\documentclass[12pt]{amsart}

\usepackage{amssymb} 
\usepackage[mathscr]{eucal} 
\usepackage[all]{xy}
\usepackage{amsmath} 

\def\cal#1{{ \mathcal {#1} }} 
\def\frak#1{{ \mathfrak {#1} }} 

\newcommand{\tagnum}[1] 
{\noindent \hbox to 22pt {\hss\bf #1.}}

\theoremstyle{plain}

\newtheorem{thm}{Theorem}[section]
\newtheorem{cor}[thm]{Corollary}
\newtheorem{lem}[thm]{Lemma}
\newtheorem{prop}[thm]{Proposition}
\newtheorem{defprop}[thm]{Definition-Proposition}

\newtheorem{defn}[thm]{Definition}

\newtheorem{rem}[thm]{Remark}

\newtheorem{example}{Example}[section]

\numberwithin{equation}{section}

\def\rar{\rightarrow}
\def\rxar#1{\xrightarrow{\kern -1pt #1}}
\def\lxar#1{\xleftarrow{#1 \kern -3pt}}

\newcommand{\CZ}{{\mathbb{C}}}
\newcommand{\KZ}{{\mathbb{K}}}

\newcommand{\Pn}[1]{{ \mathbb{P}^{#1}}\relax}

\newcommand{\trans}[1]{\,^t\!#1}

\newcommand{\QED}{\hspace*{\fill}$Q.E.D.$}

\newcommand{\ra}{\ensuremath{\rightarrow}}

\def\eea{\end{eqnarray*}}
\def\bea{\begin{eqnarray*}}

\def\C{{\mathbb{C}}}
\def\E{{\mathcal{E}}}
\def\F{{\mathcal{F}}}
\def\G{{\mathcal{G}}}
\def\H{{\mathcal{H}}}

\def\I{{\mathcal{I}}}

\def\P{{\mathcal{P}}}

\def\T{{\mathcal{T}}}
\def\Z{{\mathbb{Z}}}
\def\hol{{\mathcal{O}}}
\def\cO{{\mathcal{O}}}
\def\PP{{\mathbb{P}}}
\def\i{{\iota}}
\def\fA{{\mathfrak{A}}}
\def\fB{{\mathfrak{B}}}
\def\fM{{\mathfrak{M}}}

\DeclareMathOperator{\rk}{rk}

\newcommand{\shExt}{\cal Ext\kern 2pt}
\newcommand{\shHom}{\cal Hom\kern 2pt}
\newcommand{\Syz}{\cal Syz\kern 2pt}
\DeclareMathOperator{\Hom}{Hom}

\DeclareMathOperator{\Proj}{Proj}

\DeclareMathOperator{\coker}{coker}

\DeclareMathOperator{\pd}{pd}

\newcommand{\Proof}{{\it Proof. }}

\def\hd#1{\hbox to 19pt{\hfil$#1$\hfil}} 

\def\hdotrulefill{\cleaders\hbox
{\hbox to 2pt{\cleaders\hrule\hfill}\kern 2pt}
\hfill}

\def\vdotrule{\hbox to 0 pt{\hfill
\vbox to 10.5pt{\kern 1pt
\hbox to 0pt{\hfil \vrule height2pt \hfil} \kern 2pt
\hbox to 0pt{\hfil \vrule height2pt \hfil} \kern 2pt
\hbox to 0pt{\hfil \vrule height2pt \hfil} \kern 2pt
\hbox to 0pt{\hfil \vrule height2pt \hfil} \kern 2pt
\vss}\hfill}}

\title[A remarkable moduli space of rank 6 vector bundles]{A
remarkable moduli space of rank 6 vector bundles
related to cubic surfaces}

\author{Fabrizio Catanese}

\address{Fabrizio Catanese, Lehrstuhl Mathematik VIII,
Mathematisches Institut der
Universit{\"a}t  Bayreuth, NW II \\
D-95440 Bayreuth,  Deutschland.}

\author{Fabio Tonoli }

\address{Fabio Tonoli, Dipartimento di Matematica Univ. di Trento\\
Via Sommarive, 14 38050-Trento (Italy)}

\email{fabrizio.catanese@uni-bayreuth.de}
\email{tonoli@science.unitn.it}

\date{\today}

\begin{document}

\begin{abstract}
We study the moduli space $\fM^s(6;3,6,4)$ of simple rank 6 vector
bundles $\E$ on $\PP^3$ with Chern polynomial $1+3t+6t^2+4t^3$
and properties of these bundles, especially we prove some
partial results concerning their stability.
We first recall how these bundles are related to the construction of
sextic nodal surfaces
in $\PP^3$ having an even set of 56 nodes (cf. \cite{CaTo}).
We prove that there is an open set, corresponding to the simple
bundles with minimal cohomology, which is irreducible of dimension 19
and bimeromorphic to an open set $\fA^0$ of the G.I.T. quotient space
of the projective space
$\fB:=\{B\in \PP(U^\vee\otimes W\otimes V^\vee)\}$
of triple tensors of type $(3,3,4)$ by the natural action of
$SL(W)\times SL(U)$.

We give several constructions for these bundles,
which relate them to cubic surfaces in 3-space $\PP^3$
and to cubic surfaces in the dual space $(\PP^3)^{\vee}$.
One of these constructions,
suggested by Igor Dolgachev, generalizes to other types of tensors.

Moreover, we relate the socalled {\em cross-product involution}
for  $(3,3,4)$-tensors, introduced in \cite{CaTo}, with
the Schur quadric associated to a cubic surface in $\PP^3$ and study
further properties of this involution.
\end{abstract}

\maketitle

\section*{Introduction}
A good motivation for the study of moduli spaces of vector bundles in
$\PP^3$ comes from the classical problem concerning the geometry of
nodal surfaces $F$ in $\PP^3$, and more specifically from the study
of even sets $\Delta$ of nodes on them (Beauville has shown in \cite{Bea}
that a surface of low degree and with many nodes contains necessarily
several such even sets  $\Delta$).

In turn, as we recall in section 1 of the paper, every even set occurs
as the second degeneracy locus of a symmetric map of a vector bundle
in $\PP^3$: more precisely, the main theorem of \cite{CaCa}
asserts that, given
$\delta \in \{0,1\}$ and a
$\delta$-even set of nodes $\Delta$ on a nodal surface $F$ of
degree $d$ in $\PP^3$, there is a vector bundle $\E''$ on $\PP^3$
and a symmetric map ${\E ''}^\vee(-d-\delta) \rxar\varphi \E'' $
such that $\Delta $ is precisely the locus where the corank of
$\varphi$ equals 2, while $ F \setminus \Delta$ is the locus
of points where corank ($\varphi$) = 1.

The simplest case is the case where the vector bundle is a direct
sum of line bundles (cf. \cite{Ca}): here $\varphi$  is just
a symmetric matrix with entries homogeneous polynomials of fixed
degrees.

The problem of existence and classification of even sets of nodes
is in general based on a preliminary analysis of certain moduli spaces
of pairs $( \E'', \varphi )$ as above.  The main question being whether
a pair as above defines loci $F : = \det  (\varphi)$ and $\Delta :=
\{ x | {\rm corank} ( \varphi ) \geq  2 \}$ which have the desired
singularities.

For instance, the classification of $0$-even nodes on sextic surfaces
was achieved in \cite{CaTo}, showing the existence of sets of
cardinality 56. In this case a remarkable feature was that the
corresponding moduli space of vector bundles $\E : = \E''(3)$
   was shown
to be irreducible, yet  computer experiments showed that for a general
such bundle $\E$ the determinant
$ F : = \det \varphi$ would
yield a cubic surface $G$ counted with multiplicity
two and independent
of  the choice of $\varphi$.

This was the first motivation to try to understand the relation occurring
between such bundles and cubic surfaces. This brought to a finer
investigation of the moduli space of pairs $(\E, \varphi )$ as above,
which revealed that the latter moduli space is indeed reducible,
with a second component of codimension 7 corresponding to
cubic surfaces $G$ reducible as the union of a plane and a
smooth quadric intersecting transversally.
Several explicit constructions for the vector bundles $\E$
in question allowed finally to show that for
a general pair $ (\E, \varphi )$
in the second component   the determinant of the symmetric map $
\varphi$ yields a nodal surface $F$ with an even set $\Delta$ as wished
for.

The purpose of the present  paper is threefold: first of all we want
to  explain
the beautiful geometry relating one of our vector bundles $\E$
to a cubic surface $G$ in $\PP^3$ and to another cubic surface
$G^*$ in the dual projective space.

This relation goes back to the classical discovery that a smooth cubic
surface can be written as the determinant of a matrix of linear forms,
and also as the image of the plane $\PP^2$ through the linear system
of cubics passing through 6 given points.

Since our second purpose, following a suggestion of Igor Dolgachev,
is to show how the above correspondence can be vastly generalized (this
is done in section 2), we try to set up the classical story in a context of
modern multilinear algebra.

Let $B$ be a  general tensor of type $(3,3,4)$, more precisely let
$ B \in U^\vee\otimes W\otimes V^\vee$, where
$U,W,V$ are complex vector spaces of  respective dimensions 3,3,4.

Now, it is classical that to a  general $ 3 \times 3 \times 4$ tensor $B$
one can associate a cubic surface in $\PP^3$
by taking the determinant of the corresponding
$3\times 3 $ matrix of linear forms on $\PP^3$.
In this way we get a cubic surface $G^*$
in the dual projective space ${\PP^3}^\vee =
   \PP (V) : = \Proj (V^\vee) $,
together with two different realizations of $G^*$
as a blow up of a projective plane $\Proj (U^\vee) $
(respectively $\Proj (W)$) in a set of 6 points.
These are the points where the $3\times 4$ Hilbert--Burch
matrix of linear forms on $U$ (respectively on $W^\vee$)
drops rank by 1,
and the rational map to ${\PP^3}^\vee$  is given
by the system of cubics through the 6 points,
a system which is generated by the determinants
of the four $ 3 \times 3$ minors
of the Hilbert--Burch matrix.
One passes from one realization to the other
simply by transposing the tensor,
and we shall call this the {\bf trivial involution}
for  $3 \times 3 \times 4$ tensors.

For a $ 3 \times 3 \times 4$ tensor 
$B\in  U^\vee \otimes W\otimes V^\vee$,
besides this trivial involution, which consists in
regarding $B$ as element of $W\otimes U^\vee \otimes V^\vee$,
there exists another involution,
called the {\bf cross-product involution} (see \cite{CaTo}).
This second involution associates to a general tensor
$B\in  U^\vee \otimes W\otimes  V^\vee $ another tensor
$B '\in {U'}^\vee\otimes W\otimes V$,
where $U'$
is defined as the kernel of the map
$\Lambda^2(W^\vee)\otimes V \ra  U^\vee \otimes W^\vee $
induced by contraction with $B$.
The {\bf reversing construction}\footnote{%
We are pedantic with the order of the spaces of
a multitensor, but this is essential for a correct
correspondence between the various constructions we will develop
from a multitensor.} is then defined as 
the composistion of the cross-product involution with the trivial
involution, and associates to 
$B\in U^\vee \otimes W\otimes V^\vee$ the tensor
$B '\in W\otimes {U'}^\vee \otimes V$,

In the paper \cite{CaTo} the authors give the following direct geometric
construction of nodal sextic surfaces with an even sets of 56 nodes.

Consider the open set $\fB^*$ of $\fB=\PP(U^\vee \otimes W \otimes
V^\vee )$ given by the $(3,3,4)$-tensors $B$ whose determinant
(as $3\times3$ matrices) defines a cubic surface
$G^*\subset{\PP ^3}^\vee= \Proj(V^\vee)$.
To such a B we apply the reversing construction and we consider
the following exact sequence induced by the $(3,3,4)$-tensor
$ B '$ on $\PP (V^{\vee}) = \Proj(V)$:
$$0 \rar W^\vee \otimes \cal O (-1)  \rxar{ B '}{ U'}^\vee \otimes \cal O
\rar \cal G \rar 0.$$
Observe that if $ B '$ never drops rank by two then $\cal G$ is
an invertible sheaf on the cubic surface $G$ associated to $ B '$.

The direct construction produces  a bundle $\E$
as an extension of $6\cO$ by the sheaf
$\tau:=\cal G^{\otimes  2}(-1)$:
it turns out that the
extension $\E$ is unique (up to isomorphisms)
if the cubic surface $G$ is
smooth, and also if it is reducible as the transversal
union of a plane with a smooth quadric. 
In the latter case we obtain a sextic nodal surface $F$  as
the corank 2 degeneracy locus  of a general symmetric map
$\varphi:\ \cal E^\vee \ra \cal E$
(while, if $G$ is smooth, $F$ is
the cubic surface $G$ counted with multiplicity two)
(cf. again \cite{CaTo}).

\medskip
Concerning the cross-product involution, we show its
relation to the Schur quadric.
We fix an orientation of $W$, which allows us to identify
$\Lambda^2(W^\vee)$ with $W$.
The  Schur quadric $Q$ is defined (up to scalars)
as a generator of the kernel of the natural map
$S^2V \ra \Lambda^2 (U^\vee) \otimes \Lambda^2 W$
obtained as the composition of $S^2(B)$ with
the projection of
$$S^2(U^\vee \otimes W)=
\left(\Lambda^2 (U^\vee) \otimes \Lambda^2 W \right)
\oplus \left(S^2(U^\vee)\otimes S^2W\right)$$
onto the first factor.
Indeed, $\dim S^2V=10$, $\dim \left(\Lambda^2 (U^\vee) \otimes
\Lambda^2 W \right)=9$, and the kernel is 1-dimensional, cf. \cite[$\S 0$
and Thm 0.5]{dk}.

Classically, the Schur quadric induces an isomorphism $q$
between $\PP (V^\vee)$ and $\PP (V)$ sending
the cubic surface $G$ to the cubic surface $G^*$:
here, we consider the tensor $B$ as inducing
an injective map  $ U \ra  W \otimes V^\vee$, and
we show that the subspace
$U'$ is equal to  the image inside $W \otimes V$
of the composition of this inclusion with
$id_W \otimes q : W \otimes V^\vee \ra W\otimes V$.
Hence we obtain the tensor $ B'$ associated to the
inclusion $ U' \ra W \otimes V$.

\medskip
Up to now we have been talking about moduli spaces,
   for instance about the
moduli space of the  bundles
$\E$ which we obtain from our tensors $B$.
The trouble however is: does such a moduli space really exist?
The answer is positive because we show that the
bundles $\E$ are simple
rank 6 vector bundles with Chern polynomial
$1+3t+6t^2+4t^3$, and from the holomorphic point of view one
has a moduli space of simple vector bundles.
   We show that our
construction leads then to the realization of an open set
   of the moduli space
$\fM^s(6;3,6,4)$  of simple rank 6 vector bundles $\E$ on $\PP^3$ with
Chern polynomial
$1+3t+6t^2+4t^3$. This open set is the biholomorphic image
of an open set of the G.I.T. quotient space $\fB / SL (U) \times SL (W) $
of
$\fB:=\{B\in
\PP(U^\vee\otimes W\otimes V^\vee)\}$.
The above is a subset of the open set formed by the simple bundles with
minimal cohomology, as explained in the following:

\medskip
{\bf Main Theorem} {\em
Consider the moduli space $ \frak M^s ( 6; 3,6,4)$ of
rank 6 simple vector
bundles $\cal E$ on $\Pn 3 : = \Proj (V)$
with Chern polynomial $ 1 + 3t + 6 t^2 + 4 t^3$,
and inside it the  open set  $\ \frak A$  corresponding to the simple
bundles with minimal cohomology, i.e., those with
$$
\begin{array}[c]{llll}
(1) & H^i(\E) = 0 \ \forall i \geq 1;
\qquad & (2) & H^i(\E (-1)) = 0 \ \forall i \neq 1;\cr
(3) & H^i(\E (-2)) = 0 \ \forall i \neq 1;
\qquad & (4) &H^i(\E (-3)) = 0 \ \forall i;\cr
(5) & H^i(\E (-4)) = 0 \ \forall i.
\end{array}
$$
%

Then   $\frak A$  is irreducible of dimension 19
and it is bimeromorphic to $ \frak A^0 $,
where $\frak A^0$ is an open set of the G.I.T. quotient space of the
projective space $\frak B$ of
tensors of type $(3,4,3)$,
$\frak B : = \{B\in \PP ({U}^\vee\otimes W \otimes V^\vee)\}$
by the natural action of $SL(W) \times SL(U)$ (recall that $U,W$ are two
fixed vector spaces of dimension 3, while $ V = H^0( \PP^3, \hol (1))$.

Let  moreover $ [B] \in \frak A ^0$ be a general point:
then to $[B]$
corresponds a vector bundle $\E_B$ on $\PP^3$ 
via the kernel construction, 
and also a vector bundle $\E^*_B$ on ${\PP^3}^\vee$, 
obtained from the direct construction
applied to the tensor $B\in U^\vee \otimes W \otimes V^\vee$
(cf. definition \ref{dir_cos} applied to $B$, or equation (\ref{*})).
$\E^*_B$ is the vector bundle $\E_{B'}$, 
where $B'\in W \otimes {U'}^\vee \otimes V$
is obtained from $B$ via the reversing construction
and $[B'] \in {\frak A ^0_*}$.
}\medskip

Moreover, in section 5 we address the question of
Mumford-Takemoto (slope), respectively Gieseker (semi)stability
of the bundles $\E$. This is an interesting but delicate question,
since slope stability implies Gieseker stability, which implies slope
semistability. We can prove slope semistability, but there remains the
interesting question whether the general bundle $\E$ associated
to a tensor $B$ is Gieseker stable.

\section{Quadratic sheaves, nodal surfaces, and related vector bundles}

The study of vector bundles of  rank 6 is a slightly unusual topic of
research, in the sense that a topic of this type is usually
studied not for its own sake, but in view of applications
to other problems.
Since this is exactly the case here,   in this section
    we want to explain how we  got interested
in our class of vector bundles and in which context we encountered them.

\begin{defn}
Let $F$ be a locally Cohen--Macaulay projective scheme.
A locally Cohen--Macaulay, reflexive, coherent $\cO_F$- sheaf $\F$ is
said to be
$\delta/2$-quadratic ($\delta\in\{ 0,1 \}$)
if there exists a symmetric isomorphism
$\sigma:\ \F(\delta) \ra \H om_{\cO_F}(\F,\cO_F)$
(symmetric means that the associated bilinear map
$ \F \times \F \ra \cO_F(- \delta)$ is symmetric).
\end{defn}

Let us now suppose that $F$ is a hypersurface of degree $d$ in a

projective space $\PP$.
If $\F$ is a coherent $\cO_F$- sheaf, we  have a natural isomorphism:
$$\H om_{\cO_F}(\F,\cO_F(-\delta)) \cong
\E xt ^1_{\cO_{\PP}}(\F,\cO_{\PP}(-d-\delta)).$$
Indeed, let $f=0$ be the equation defining $F$ and consider the exact
sequence
$$0\ra \cO_{\PP}(-d-\delta) \rxar{\cdot f}
\cO_{\PP}(-\delta) \ra \cO_F(-\delta)\ra0.$$
By applying the functor $\H om_{\cO_{\PP}}(\F,-)$,
we obtain
$$0 \ra \H om_{\cO_F}(\F,\cO_F(-\delta)) \ra
\E xt ^1_{\cO_{\PP}}(\F,\cO_{\PP}(-d-\delta)) \rxar{\cdot f}
\E xt ^1_{\cO_{\PP}}(\F,\cO_{\PP}(-\delta)),$$
where the last map is zero, since it is induced by  multiplication
by $f$ as a morphism between $\cO_F$- sheaves of modules.

Therefore, for a quadratic sheaf $\F$ defined on $F$,
we have
$$\F \cong \H om_{\cO_F}(\F,\cO_F(-\delta))\cong
\E xt ^1_{\cO_{\PP}}(\F,\cO_{\PP}(-d-\delta)).$$

On the other hand,
if $\E''$ is a vector bundle on $\PP$ and
$\varphi:\ {\E''}^\vee(-d-\delta) \ra {\E''}$ is a symmetric morphism,
we define $F$ as the locus where $\rk(\varphi)\leq\rk{\E''}-1$
and set $\F:=\coker\varphi$.
We then have the exact sequence
$$0 \ra {\E''}^\vee(-d-\delta) \ra {\E''} \ra \F \ra 0$$
and, by applying $\H om_{\cO_{\PP}}(-,\cO(-d-\delta))$,
$$0 \ra {\E''}^\vee(-d-\delta) \ra {\E''} \ra
\E xt ^1_{\cO_{\PP}}(\F,\cO_{\PP}(-d-\delta)) \ra 0.$$
Thus $\F$ is naturally isomorphic to
$\E xt ^1_{\cO_{\PP}}(\F,\cO_{\PP}(-d-\delta))$.
By identifying $\E xt ^1_{\cO_{\PP}}(\F,\cO_{\PP}(-d-\delta))$
with $\H om_{\cO_F}(\F,\cO_F(-\delta))$ via the
natural isomorphism  described above,
we  get a symmetric isomorphism
$\F \cong \Hom(\F,\cO_F(-\delta))$,
and we finally conclude that $\F$ is a quadratic sheaf.
Thus a symmetric morphism
$\varphi:\ {\E''}^\vee(-d-\delta) \ra {\E''}$ induces a quadratic sheaf
with support on the hypersurface $F:=\{\det \varphi=0\}$.

Assume now that the generic rank of $\F$ is $1$, and choose a section
$\beta \in H^0 (\hol_F (\delta))$:
then the bilinear pairing  $ \F \times \F \ra \cO_F(- \delta)$ composed
with multiplication by $\beta$ yields a commutative ring structure on
the module $\cO_F \oplus \F$.

We can then
    consider the scheme
$\widetilde{F}:=Spec(\cO_F\oplus  \F)$,
which yields a 2:1 covering $\pi:\ \widetilde{F}
\ra F$,
\'etale over the complement of $\Delta \cup \{ \beta = 0\}$, where
$\Delta$ is the locus where
$\rk(\varphi)\leq\rk\E-2$.

In this way, a subscheme $\Delta$ of  a locally
Cohen--Macaulay projective scheme $F$
is said to be bundle--symmetric if there exist a bundle $\E$ and a
symmetric morphism $\varphi$ such that $F$
is the corank 1 locus of $\varphi$
and $\Delta$ is the corank 2 locus of $\varphi$, as in the above
setting.

\begin{defn}
Let $F$ be a nodal surface in $\PP^3$ and
$\Delta$ be a set of nodes of $F$.
Consider the resolution $\widetilde F$ of $F$ along the singularities
in $\Delta$  and the corresponding $(-2)$-curves $A_i$.
A set of nodes $\Delta\subset F$ is called $\delta/2$-even if
the corresponding divisor $\sum A_i+\delta H$
is 2-divisible in the Picard group of $\widetilde F$,
where $H$ is the divisor class corresponding to a hyperplane section
of $F$ in $\PP^3$ (cf. \cite{Ca}).
\end{defn}

This property amounts  to the existence of a $2:1$ covering $\widetilde
S$ of $\widetilde F$ ramified along the divisor $\sum A_i+\delta H$,
or, equivalently, to the existence of a $2:1$ covering $S$ of $F$
ramified on $\Delta$ (and on a hyperplane section if $\delta = 1$)
($S$ is obtained by blowing down the $(-1)$ rational
curves on $\widetilde S$ which are the inverse images of the $ A_i$'s).

\medskip
Clearly, a bundle--symmetric set of nodes is a $\delta/2$-even set,
but also the converse holds, as it was shown in \cite{CaCa}:

\begin{thm}\cite[theorem 0.3 and corollary 0.4]{CaCa}\label{caca}
    Let $F\subset\PP^3$ be a nodal surface of degree $d$.
Then every $\delta/2$-even set of nodes $\Delta$ on $F$, $\delta=0,1$,
is the degeneracy locus of a symmetric map of locally free
$\cO_{\PP^3}$-sheaves 
${\E''}^\vee (-d-\delta) \rxar{\varphi} \E''$
(\/i.e., $F$ is the locus where $\rk(\varphi)\leq\rk\E''-1$,
$\Delta$ is the locus where $\rk(\varphi)=\rk\E''-2$\/).
Moreover, if $p:\ S\ra F$ is a $2:1$ covering associated to
$\Delta$ and $\F$ is defined as the anti-invariant part
of $p_*(\cO_S)=\cO_F\oplus\F$, 
then $\F$ fits into the exact sequence
$$
0\ra {\E''}^\vee(-d-\delta) \rxar\varphi {\E''} \ra \F \ra 0.
$$
\end{thm}

The authors of \cite{CaCa} also describe how to construct the bundle $\E''$.
Recall that the first syzygy bundle $Syz^1(M)$ associated to
a graded module $M$ is obtained as follows.
Take a  graded free resolution of the module $M$:
$$ 0 \rar \P^4 \rar \dots \rar P^1 \rxar{\alpha_1} P^0
\rxar{\alpha_0} M \ra 0.$$
Then the homomorphism $\alpha_1: P^1 \rar P^0$ induces a corresponding
homomorphism $(\alpha_1)^\sim$ between the (Serre-) associated sheaves
$({P^1})^\sim$ and
$({P^0})^\sim$, and the first syzygy
bundle of $M$ is defined as $Syz^1(M) : = {\rm Ker} ( \alpha_1^\sim )$.

In \cite{CaCa} it is shown that, up to direct sum of line bundles,
$\E''$  is the first syzygy bundle of the module
$$
M=U\oplus \left(\oplus_{m>(d-4+\delta)/2} H^1(F,\F(m))\right),
$$
where, if $(d+\delta)$ is even,
$U$ is any maximal isotropic subspace of $H^1(F,\F((d-4+\delta)/2))$
with respect to the
non-degenerate alternating form on $H^1(F,\F((d-4+\delta)/2))$
induced by Serre's duality,
and, if $(d+\delta)$ is odd, $U$ is zero
(cf. \cite{CaCa} for more details).

\bigskip
Even sets of nodes are classified and explicitly described for surfaces
of degree up to 5.
In \cite{CaTo}, we  studied the case of  even (i.e., 0-even)
sets on sextic
surfaces.

Particularly interesting is the case of even sets of cardinality 56:
this is the first case where the module $M$ is relatively big.
   Concerning the possible dimensions of the
various graded pieces   of $M$, we showed that only   two  cases can
occur: either $h^1(F,\F(1))=h^1(F,\F(2))=3$ or $h^1(F,\F(1)) = 3,
h^1(F,\F(2))=4$. Then we studied the former case.

In the former case $U$ is thus an (isotropic) 3-dimensional vector space in
$H^1(F,\F(1))$ and, if we denote by
   $W$ the 3-dimensional vector space $W:=H^1(F,\F(2))$
and set $V : =H^0(\PP^3,\cO_{\PP^3}(1))$,
we have that $M$ is completely determined by the multiplication
tensor
\begin{equation}\label{tritensorB}
B\in U^\vee\otimes W \otimes V^\vee .
\end{equation}

We now describe Beilinson's table for $\E''$, up to direct sums of line
bundles.
Since $H^0(\PP^3,\F(2))=0$
(as shown in prop. 2.4 of \cite{CaTo}, compare  also Beilinson's table for $\F$
given in section 3 of \cite{CaTo}),
one has $H^0(\PP^3,\E''(2))=0$.
Since, up to direct sums of lines bundles,
$\E''$ is the first syzygy bundle of $M$, we have that
$H^1_*(\PP^3,\E)\cong M$ and all the other intermediate cohomology
modules of $\E''$ are zero.
It follows that  Beilinson's table $h^i(\E''(j))$ for $\E''$,
up to direct sums of line bundles, is
$$
\vbox{\offinterlineskip
\halign{
\hbox to 10pt{\hfil}#
&\vrule height10.5pt depth 5.5pt#
&\hbox to 35pt{\hfil$#$\hfil}
&\vrule#
&\hbox to 35pt{\hfil$#$\hfil}
&\vrule#
&\hbox to 35pt{\hfil$#$\hfil}
&\vrule#
&\hbox to 40pt{\hfil$#$\hfil}
&\vrule#
&\hbox to 40pt{\hfil$#$\hfil}
&\vrule#
&#
\cr
\omit&\omit&&\omit
\hbox to 0pt{\hbox to 0pt{\hss$\big\uparrow$\kern
-0.4pt\hss}\raise5pt\hbox to 0pt{\ $i$\hss}}\cr
&\multispan{2}\hrulefill&\multispan{8}\hrulefill\cr
&& 0 && 0 && 0 && 0 && 0 && \cr
&\multispan{2}\hrulefill&\multispan{8}\hrulefill\cr
&& 0 &&  0 && 0 && 0 && 0 && \cr
&\multispan{2}\hrulefill&\multispan{8}\hrulefill\cr
&& 0 &&  0 && 3 && 3 && 0 && \cr
&\multispan{2}\hrulefill&\multispan{8}\hrulefill\cr
&& 0 && 0 && 0 && 0 && * &&\cr
\multispan{12}\hrulefill&\kern -2pt
\hbox{\kern .4pt \vbox to 0pt{\kern 1.1pt\vss\hbox to
0pt{$\longrightarrow$ \hss}\vss}
                 \vbox to 0pt{\kern 3pt \hbox to 0pt{\ $j$ \hss}\vss}}\cr
}\kern 6pt}
$$
If we denote from now on by $\E$ the previous $\E(3)$,
theorem \ref{caca} gives the exact sequence
$$ 0 \ra \E^\vee \rxar\varphi \E \ra \F(3) \ra 0,$$
and, by Beilinson's theorem and the above cohomology table,
$\E(-1)$ is obtained by adding a direct sum of line
bundles to the bundle
$$
\cal E'(-1):=\ker\left(U\otimes\Omega^1(1)\cong 3\Omega^1(1)
\rxar{\beta(-1)}
W\otimes\cal O\cong 3\cal O\right).
$$
Moreover, since  Beilinson's complex has no cohomology in degree
$\neq 0$, the above map is surjective: hence
$\cal E'$ is a vector bundle with $rk(\cal E')=6$.

Recall  the Euler sequence
\begin{equation}\label{euler}
0\rar \Omega^1(1)\rar V\otimes\cal O\cong 4\cal O\rxar{\epsilon}\cal
O(1)\rar 0,
\end{equation}
where $V:=H^0(\Pn 3,\cal O(1))$
is the space of linear forms on $\Pn 3$,
and suppose that the map $H^0(\beta)$, induced in cohomology by
$\beta$, is also surjective (as it happens for a general $\beta$):
then $h^0(\PP^3,\E')=3h^1(\Omega^1(2))-12=6$.

Since $h^0(\PP^3,\F)=6$
(cf. again Beilinson's table for $\F$ in section 3 of \cite{CaTo}),
one can make the following generality assumption:

\medskip\noindent
{\bf  FIRST ASSUMPTION:  $\cal F$ is generated in degree $3$}
and the linear map $H^0(\cal E') \ra H^0(\cal F(3)) $
is an isomorphism.
\medskip

Under the above assumption  the following holds:

\begin{prop}\label{1as}\cite[Prop 3.3]{CaTo}
Notation being as above, if
the  first assumption holds true then $\cal E=\cal E'$ or, equivalently,
{\bf  rank   $(\cal E)$  = 6}. More precisely,
there exists a homomorphism $\beta :  U\otimes\Omega^1(2) \cong
3\Omega^1(2) \ra  W\otimes\cal O(1) \cong 3\cal O(1) $
with $\cal E =
ker
\beta$ and  such that we have an exact sequence
\begin{equation}
0\rar \cal E \rar U\otimes\Omega^1(2) \cong
3\Omega^1(2)  \rxar\beta W\otimes\cal O(1) \cong 3\cal O(1) \rar 0.
\end{equation}

Conversely, if $\cal E$ is obtained in this way it is a rank 6 bundle
and, if the map $H^0(\beta)$ is surjective, 
it has an intermediate cohomology module $M:=H^1_*(\E)$ 
having the required Hilbert function of type $(3,3)$.
\end{prop}

We are now able to explain why we got interested
in vector bundles $\E$ of rank 6 with
Chern polynomial $c(\E) (t) = 1 + 3 t + 6 t^2 + 4 t^3$:

\begin{lem}\label{E=beta}
Let $\E$ be given as kernel of a surjective homomorphism
$\beta :  U\otimes\Omega^1(2) \cong
3\Omega^1(2) \ra  W\otimes\cal O(1) \cong 3\cal O(1) $,
where $U$ and $W$ are 3-dimensional vector spaces:
\begin{equation}\label{E=ker}
0\rar \cal E \rar U\otimes\Omega^1(2) \cong
3\Omega^1(2)  \rxar\beta W\otimes\cal O(1) \cong 3\cal O(1) \rar 0.
\end{equation}
Then $\E$ is a rank 6 bundle with total Chern class
\begin{equation}\label{c(t)}
         c(\E) (t) = 1 + 3 t + 6 t^2 + 4 t^3,
\end{equation}
and $H^0(\E^\vee)=0$.

Moreover, if the map $H^0(\beta)$ is surjective, 
$\E$ is a bundle with minimal cohomology, more precisely, it satisfies:
$$
\begin{array}[c]{llll}
(1) & H^i(\E) = 0 \ \forall i \geq 1;
\qquad & (2) & H^i(\E (-1)) = 0 \ \forall i \neq 1;\cr
(3) & H^i(\E (-2)) = 0 \ \forall i \neq 1;
\qquad & (4) &H^i(\E (-3)) = 0 \ \forall i;\cr
(5) & H^i(\E (-4)) = 0 \ \forall i.
\end{array}
$$
In this case, in particular, $h^0(\E)=6$ and the unique nonzero
intermediate cohomology groups of $\E$ are $U=H^1(\PP^3,\E(-2))$ and
$W=H^1(\PP^3,\E(-1))$.

Conversely, a bundle $\E$ with total Chern class as above and with
minimal cohomology as above admits  a presentation  of type (\ref{E=ker}),
where
$U=H^1(\PP^3,\E(-2))$ and $W=H^1(\PP^3,\E(-1))$ are
3-dimensional vector spaces and $H^0(\beta)$ is surjective.
\end{lem}

\proof
If $\E$ is given as in (\ref{E=ker}), then it is a rank 6 bundle
and its Chern polynomial is
$c(\cal E)=c(\Omega^1(2))^3 c(\cal O(1))^{-3}=
(c(\cal O(1))^4 c(\cal O(2))^{-1})^3c(\cal O(1))^{-3}=
(1+t)^9(1+2t)^{-3}=(1+9t+36t^2+84t^3)(1-6t+24t^2-80t^3)=
1+3t+6t^2+4t^3$.

Dualizing the sequence
$0\rar \cal E \rar 3\Omega^1(2) \rxar\beta 3\cal O(1) \rar 0$
yields
$$  0\rar 3\cal O(-1) \rar 3 T(-2) \rar {\cal E}^\lor \rar 0.$$
Thus $h^0( {\cal E}^\lor)=0$.

\medskip
Let us verify that a bundle given as in
(\ref{E=ker}) satisfies properties (2) and (3).
The exact cohomology sequences of the twists of (\ref{E=ker}) give: 
$H^1 ( \cal E (-2))\cong 3 H^1 (\Omega^1)$, 
$H^1 ( \cal E (-1)) \cong 3 H^0 (\hol)$.
This also shows that all the other intermediate cohomology modules
of the above twists of $\E$ are zero.
Considering the Euler sequence,
it is also clear that all the negative twists of $\E$
have no global sections, and that $ \cal E (-2)$ and  $\cal E (-1)$ 
have vanishing third cohomology groups.
Therefore $ \cal E (-2)$ and  $\cal E (-1)$, 
if  $\E$ is given as in (\ref{E=ker}),
have only first cohomology group, and of dimension 3.

The exact cohomology sequence of (\ref{E=ker}) gives 
$H^1( \E)=\coker(H^0(\beta))$ and $H^2(\E)=0$.
If $H^0(\beta)$ is assumed to be surjective, then also property (1) is 
satisfied and $h^0(\E)=\chi(\E)=6$ is determined by the Riemann-Roch theorem.

Moreover, if $H^0(\beta)$ is surjective, then also 
$H^0(\beta(k))$ is surjective for positive twists $k$.
Since $\Omega^1(k)$ has no global section for $k\leq 0$, 
it follows easily that the intermediate cohomology group
$H^1_*(\E)$ has nonzero degree parts only in degree -2 and -1
and that $H^2_*(\E)=0$. 

It remains to show the vanishing of the groups $H^3(\E(k))$,
for $-4\leq k\leq 0$.
This follows from the vanishing of $H^3(\Omega^1(k))$, for $-2\leq
k\leq 2$, which is a straighforward computation:
$H^3(\Omega^1(k))\cong H^0(\T^1(-k-4))=0$ for $-2\leq k\leq 2$.

We now prove the converse.
Assume that we have a vector bundle with such a Chern polynomial
and  minimal cohomology as described above.
Then the Euler characteristics of $\E$ (or its twists) are the
same as the Euler characteristic of a bundle (or its twists)
given as in (\ref{E=ker}). 
Therefore the first cohomology groups
$U=H^1 ( \cal E (-2))$ and $W=H^1 ( \cal E (-1))$ are
both 3-dimensional vector spaces.

By appling Beilinson's theorem to $\E(-1)$,
it follows that $\E$ 
fits in an exact sequence as in (\ref{E=ker}).
Condition (1) implies that $H^0(\beta)$ is surjective.
\qed

%

\section{Sheaves associated to tensors}
Let $W_1,\dots,W_r$ be vector spaces of respective dimensions
$\dim_\KZ W_j =d_j+1$, where $\KZ$ is $\CZ$ or any algebraically
closed field.
Assume that we have a tensor
\begin{equation}
B\in W_1\otimes \dots \otimes W_r.
\end{equation}

To $B$ one can associate a collection of subschemes
of products of projective spaces, namely
\begin{defn}
     Let $1\leq i_1< \dots < i_h\leq r$ be a strictly increasing
$h$-tuple of indexes between $1$ and $r$. Then
we define
$$
\Gamma_{i_1,\dots,i_h}(B):=\left\{
u=u_{i_1}\otimes\dots\otimes u_{i_h}\mid
u_{i_j}\in W_{i_j}^\vee,\ B\neg u=0
\right\};
$$
where $B\neg u$ denotes the contraction of the tensor $B$
with $u$.
\end{defn}

\begin{rem}
$1)$  Without loss of generality, we may assume
$i_1=1,\dots,i_h=h$, and $r=h$ or $r=h+1$
(it suffices to replace the vector spaces
$W_{h+1},\dots,W_{r}$ with their tensor product
$W_{h+1}\otimes\dots\otimes W_{r}$).

$2)$ Observe that
$\Gamma_{i_1,\dots,i_h}(B)
\subset\PP(W_{i_1}^\vee)\times \dots\times\PP(W_{i_h}^\vee)$
is the intersection of hypersurfaces of multidegree
$(1,\dots,1)$.

$3)$ Is $h=r$, then we have a single hypersurface of
multidegree $(1,\dots,1)$.
Otherwise, we shall make the following assumption of generality:
$\Gamma_{1,\dots,r-1}(B)$ is the complete intersection of
$d_r+1$ hypersurfaces of multidegree $(1,\dots,1)$.

$4)$ The case $r=1$ is empty, while the case $r=2$ is not very
interesting, since we get corresponding linear maps
$A:\ W_1^\vee \ra W_2$ and $\trans A:\ W_2^\vee \ra W_1$
and loci $\PP(\ker A)$, $\PP(\ker \trans A)$,
$\{(x,y)\in \PP (W_1^\vee ) \times \PP ( W_2^\vee )  \mid\
\langle y,A x\rangle =0\}$

$5)$ Observe finally that it suffices to treat the case
$r=h+1$.
In fact, if $r=h$, we take $W_{r+1}=\CZ$, whereas the case
$r>h+1$ can be reduced to the case $r=h+1$, as observed in
part (1).
\end{rem}

We now fix a tensor $B\in W_1\otimes \dots W_{h+1}$ as above.
In order to study the sheaves associated to the tensor $B$,
the following assumption is fundamental.

{\bf Main Assumption:}
{\em consider
$\PP:=\PP(W_{1}^\vee)\times \dots\times\PP(W_{h}^\vee)$
and the variety $\Gamma:=\Gamma_{1,\dots,h}\subset \PP$.
We assume that
$\Gamma:=\Gamma_{1,\dots,h}\subset \PP$
is a complete intersection of $d_{h+1}+1$ hypersurfaces
of multidegree $(1,\dots,1)$.

We further assume that $\Gamma\neq \emptyset$
(under the above assumption,
this happens if and only if $d_1+\dots+d_h\geq d_{h+1}+1$).}

If the main assumption holds we have then a Koszul exact sequence
{\small
\begin{equation}
     \label{stella}
     \dots
\ra \wedge^2 W^\vee\otimes \cO_\PP(-2,\dots,-2)
\ra W^\vee \otimes \cO_\PP(-1,\dots,-1)
\ra \cO_\PP \ra \cO_\Gamma \ra 0,
\end{equation}
}
where
$$
\begin{cases}
     W : =W_{h+1}\\
     \PP : =\PP(W_{1}^\vee)\times \dots\times\PP(W_{h}^\vee)\\
     \PP(W_{i}^\vee):=\Proj(Sym(W_i))
\end{cases};
$$
or, equivalently,
{\small
\begin{equation}
     \label{diesis}
     \dots
\ra \wedge^2 W^\vee\otimes \cO_\PP(-2,\dots,-2)
\ra W^\vee \otimes \cO_\PP(-1,\dots,-1)
\ra \I_\Gamma \ra 0,
\end{equation}
}

\begin{defn}
Assume now that $s<h$:
consider
$\PP':=\PP(W_{1}^\vee)\times \dots\times\PP(W_{s}^\vee)$,
$\PP''=\PP(W_{s+1}^\vee)\times \dots\times\PP(W_{h}^\vee)$,
and let $p:\ \PP \ra \PP'$ be the projection of the product
$\PP=\PP'\times\PP'' \ra \PP'$ onto the first factor.

For $\overline t=(t_{s+1},\dots,t_{h})$, we define
$\G_{\overline t}:=
p_*\cO_\Gamma(0,\dots,0,t_{s+1},\dots,t_{h})$.
\end{defn}

We aim at giving a resolution of $\G_{\overline t}$.
The above situation is quite general, but in any case
the object $\cO_\Gamma(0,\dots,0,t_{s+1},\dots,t_{h})$
can be replaced, as an object in the derived category
of coherent sheaves on $\PP$, by the twisting
by $\cO_\PP(0,\dots,0,t_{s+1},\dots,t_h)$ of the resolution
(\ref{stella}).
By applying $p_*$ to the exact sequence obtained in this
way,  we get a spectral sequence converging to
$R^k p_* \cO_\Gamma(0,\dots,0,t_{s+1},\dots,t_h)$
( if $h = s+1$ we get a complex,
as in Beilinson's theorem (cf. \cite{Bei}),
whose $k$-th cohomology group is
$R^k p_* \cO_\Gamma(0,\dots,0,t_{s+1},\dots,t_h)$).

The advantage of using the twisted Koszul complex
(\ref{stella}) is that a line bundle
$\cO(a_{s+1},\dots,a_h)$ on $\PP''$ is an external tensor
product $\cO(a_{s+1})\boxtimes \dots \boxtimes\cO(a_h)$,
hence its total cohomology
$H^*(\PP'',\cO(a_{s+1},\dots,a_h))$ is the tensor
product $H^*(\PP(W_{s+1}^\vee),\cO(a_{s+1}))\otimes\dots
\otimes H^*(\PP(W_{h}^\vee),\cO(a_{h}))$.

On the other hand,
$H^*(\PP^d,\cO_{\PP^d}(a))$ contains at most one term:
$H^0(\PP^d,\cO_{\PP^d}(a))$ if $a\geq 0$,
$H^d(\PP^d,\cO_{\PP^d}(a))$ if $a\leq -d-1$,
none if $-d\leq a\leq -1$.
Whence, fixed i,
$R^j p_* \cO_\PP(0,\dots,0,t_{s+1}-i,\dots,t_h-i)=0$
with only one possible exception $j$.

We thus obtain the following proposition.

\begin{prop}
There is a spectral sequence with $E_1$ term
$E_1 (-i,j)$ given by
\begin{multline}\label{E_1}
R^j p_* \left(
\wedge ^i W^\vee \otimes
\cO_\PP(-i,\dots,-i,t_{s+1}-i,\dots,t_h-i)\right)=\\
=\wedge^i W^\vee
\otimes H^j(\cO_{\PP''}(t_{s+1}-i,\dots,t_h-i))
\otimes \cO_{\PP'}(-i,\dots,-i)
\end{multline}
which converges to the direct image sheaves
$R^k p_* \cO_\Gamma(0,\dots,0,t_{s+1},\dots,t_h)$.

\end{prop}

\proof
This is a standard spectral sequence argument, compare pages 149-150
of \cite{weibel}.
Consider the complex given by (\ref{stella})
(without the last term at the right)
and tensor it by $\cO_\PP(0,\dots,0,t_{s+1},\dots,t_h)$.
The sequence obtained, say $C^\cdot$, remains exact,
with the exception of the rightmost term,
where the cohomology group is $\cO_\Gamma(0,\dots,0,t_{s+1},\dots,t_h)$.

If we take an injective resolution of the complex and
apply the  functor  $p_*$, we obtain a double complex with two
associated spectral sequences.
The horizontal spectral sequence degenerates at the $E_1$ term,
and yields the direct image sheaves $ E^{hor}_1 (0,k) = R^k p_*
\cO_\Gamma(0,\dots,0,t_s,\dots,t_h)$.

The vertical spectral sequence, instead, yields an $E_1$ term
of the form $E^{vert}_1 (-i,j) =   R^j p_* \left(
\wedge ^i W^\vee \otimes
\cO_\PP(-i,\dots,-i,t_{s+1}-i,\dots,t_h-i)\right).$
This is precisely the spectral sequence which we choose, and which
converges to 
$ R^k p_*\cO_\Gamma(0,\dots,0,t_{s+1},\dots,t_h)$ as claimed.

\qed

\medskip
Consider the differential $d_1$ at the $E_1$ level of the vertical
spectral sequence:
it is a horizontal differential given at the place $(-i,j)$  by
a map
\begin{multline}
\wedge^i W^\vee
\otimes H^j(\cO_{\PP''}(t_{s+1}-i,\dots,t_h-i))
\otimes \cO_{\PP'}(-i,\dots,-i)
\ra\\
\wedge^{i-1} W^\vee
\otimes H^j(\cO_{\PP''}(t_{s+1}-i+1,\dots,t_h-i+1))
\otimes \cO_{\PP'}(-i+1,\dots,-i+1)
\end{multline}
induced by \ref{stella}.

\smallskip
By the above discussion on the cohomology groups
$H^*(\PP'',\cO(a_{s+1},\dots,a_h))$,
first of all
it follows that the term $E_1 (-i,j)$ is nonzero, for fixed i,
only for at most one value of $j$.

More precisely, if $ H^j(\cO_{\PP''}(t_{s+1}-i,\dots,t_h-i)) \neq 0$,
then there is an expression $ j = j_{s+1} + \dots +j_h $
such that the above group is an external tensor product
of cohomology groups $ H^{j_c} (\PP^{d_c}, \hol (t_c-i))$.
Since each term of the
external tensor product must be nonzero,
it follows that $j_c = 0$ or$j_c = d_c$ and that
$t_c -i \geq 0$  if $j_c = 0$, else
$t_c -i \leq - d_c -1$.

\smallskip
Moreover, we conclude also that 
$H^{j-p}(\cO_{\PP''}(t_{s+1}-i+p+1,\dots,t_h-i+p+1)) = 0$
unless there is are some $ j_c = d_c $ such that $ H^{d_c} (\PP^{d_c}, \hol
(t_c-i)) \neq 0$ and $ H^{0} (\PP^{d_c}, \hol (t_c-i+ p+ 1))
\neq 0$: this is only possible if $t_c -i \leq - d_c -1$ and $t_c-i+
p+ 1 \geq 0$.
This implies $ - d_c -1 \geq t_c - i \geq -p -1 $, in particular, $ p
\geq d_c$.

\medskip
We want now to consider an easier situation, first of all we would like to have
$$
R^jp_*(\cO_\Gamma(0,\dots,0,t_{s+1},\dots,t_h)) = 0 \text{ for } j\geq 1,
$$
so that the given spectral sequence converges then to
$R^jp_*(\cO_\Gamma(0,\dots,0,t_{s+1},\dots,t_h))$.

To achieve this property, we assume $s=h-1$.

\begin{lem}\label{t>0}
If $s=h-1$, $R^jp_*(\cO_\Gamma(0,\dots,0,t))$
for $j\geq 1$, assuming that $t\geq -1$.
\end{lem}

\proof
This follows from the base change theorem, since the fibres of
$\Gamma \ra \PP'=\PP^{d_1}\times\dots\times\PP^{d_{h-1}}$
are linear subspaces of $\PP^{d_h}$,
and since $H^j(\cO_{\PP^d}(t))=0$
for any $d$, $t\geq -1$, and $j\geq 1$.

\qed

\begin{cor}\label{cor1}
If $s=r-2$ and $h=r-1$, then for $t\geq -1$
$\G_t=p_* \cO_\Gamma(0,\dots,0,t)$ has a resolution
given by an exact complex of vector bundles on
$\PP'$ whose $k$-th term is
\begin{equation}
\bigoplus_{j-i=k} \wedge^i W^\vee
\otimes H^j(\cO_{\PP^{d_h}}(t-i))
\otimes \cO_{\PP'}(-i,\dots,-i).
\end{equation}
\end{cor}

\proof
In this case the spectral sequence degenerates at the $E_2$ level
if all the nonzero terms of $E_1$ occur only for $ j=0$.

If there is a nonzero term with $ j=d_h$, then, as we observed,
for $ p > 0$ we have 
$H^{j-p}(\cO_{\PP^{d_h}}(t-i+ p + 1)) \neq 0$ iff $ p = d_h$
and $ t - i = - d_h - 1$.

In terms of differentials of the spectral sequence,
this implies that on this nonzero term 
$d_1= d_2= d_{p-1} = 0$,
and then also in the corresponding place 
$d_{p+1} =  \dots = 0$.
The result now follows.

\qed

\bigskip
We want  now to restrict ourselves to the case where we obtain sheaves on
projective spaces,
i.e., we restrict to the case $r=3$ of tritensors.

We have
$\PP:=\PP^{d_1}\times \PP^{d_2}$,
$\PP':= \PP^{d_1}$ and $p:\ \PP^{d_1}\times \PP^{d_2} \ra \PP^{d_1}$.
Recall that under the main assumption that $\Gamma\subset \PP$
be a complete intersection we have:
\begin{equation}\label{Gammanotempty}
\Gamma\neq \emptyset \iff d_1+d_2\geq d_3+1.
\end{equation}

By applying the above corollary (\ref{cor1}) we get:
\begin{cor}\label{cor2}
Suppose that $\Gamma$ is not empty, c.f. (\ref{Gammanotempty}).
Assume that $t-d_3-1\geq -d_2$ ( i.e.,  $t\geq d_3+1-d_2$)
and $ t \geq -1$.

Then $\G_t$ has a resolution given by an
Eagon--Northcott type complex:
\begin{equation}
0\ra
\begin{matrix}
\wedge^{d_3+1}W_3^\vee \\
\otimes \\ S^{t-d_3-1}W_2 \\
\otimes \\ \cO_{\PP'}(-d_3-1)
\end{matrix}
\ra \dots \ra
\begin{matrix}
     \wedge^2W_3^\vee\\
\otimes \\ S^{t-2}W_2\\
\otimes \\ \cO_{\PP'}(-2)
\end{matrix}
\ra
\begin{matrix}
     W_3^\vee\\
\otimes \\ S^{t-1}W_2\\
\otimes \\ \cO_{\PP'}(-1)
\end{matrix}
\ra
\begin{matrix}
S^{t}W_2\\ \otimes \\ \cO_{\PP'}
\end{matrix}
\ra \G_t \ra 0.
\end{equation}
\end{cor}

\begin{rem}
$1)$ Note that $t\geq d_1$ suffices in the above corollary.

$2)$ Observe that $S^t(W_2)=0$ for $t<0$, and $S^0(W_2)=\CZ$.
We may then assume $t\geq0$.
For $t=0$ we need $d_2\geq d_3+1$, and then we get
$\G_0\cong \cO_{\PP'}$.
\end{rem}

\begin{example}
     For $t=1$ we need $d_2\geq d_3$, and then we get the ``standard''
matricial resolution
\begin{equation}
     0\ra W_3^\vee\otimes\cO_{\PP'}(-1) \ra W_2\otimes \cO_{\PP'} \ra
\G_1 \ra 0,
\end{equation}
where $B\in \Hom(W_3^\vee\otimes\cO_{\PP'}(-1), W_2\otimes \cO_{\PP'})
\cong W_3\otimes W_1 \otimes W_2$.
\end{example}

\begin{example}
For $t=2$, in the case where $d_2\geq d_3-1$ we get:
\begin{equation}
     0 \ra \wedge^2 W_3^\vee\otimes\cO_{\PP'}(-2) \ra
W_3^\vee\otimes W_2^\vee\otimes \cO_{\PP'}(-1) \ra
S^2W_2\otimes \cO_{\PP'}\ra \G_2 \ra 0.
\end{equation}
\end{example}

This case is the one we are particularly interested in,
because if we set
$\widehat\E:=p_*\I_\Gamma(0,2)$, then we obtain
\begin{gather}
0 \ra \wedge^2W_3^\vee\otimes\cO_{\PP'}(-2) \ra
W_3^\vee\otimes W_2^\vee\otimes \cO_{\PP'}(-1) \ra
\widehat\E \ra 0\label{A2}\\
0 \ra \widehat\E \ra S^2W_2\otimes \cO_{\PP'}\ra \G_2 \ra 0.\label{A1}
\end{gather}

\begin{rem}
     The sheaves $\G_t$ are supported on
\begin{equation*}
\begin{split}
p(\Gamma)&=\{u_1\in W_1^\vee \mid \exists u_2\in W_2^\vee
\text{ s.t. } B\neg(u_1\otimes u_2)=0\}\\
&=\{u_1\in W_1^\vee \mid B\neg u_1 \text{ has a nontrivial kernel}\}\\
&=\{u_1\in W_1^\vee \mid \rk(B\neg u_1)\leq d_2\}.
\end{split}
\end{equation*}
In particular, the expected codimension of $p(\Gamma)$ equals
$d_3-d_2+1=d_1-\dim\Gamma$.
\end{rem}

\begin{rem}
If $p(\Gamma)$ has codimension $\geq 1$,
then $\widehat\E$ is a vector bundle
if and only if $p(\Gamma) $ is a hypersurface
and $\G_2$ is Cohen-Macaulay.
\end{rem}

\proof
Notice that (\ref{A1}) implies that, if
$\widehat\E$ is locally free, then $\G_2$
has local projective dimension at most 1 (over
the local ring $\cO_{\PP'}$). Whence
the codimension of $p(\Gamma)$ (the support of $\G_2$)
is at most 1.

Thus, if   $p(\Gamma) \neq \PP  '$,  $p(\Gamma) $ is a hypersurface.
Conversely, if $p(\Gamma) $ is a hypersurface, then
$\G_2$ is Cohen-Macaulay iff it has projective dimension 1.

We dualize the exact sequence (\ref{A1}), obtaining:
$$0 \ra S^2W_2\otimes \cO_{\PP'} \ra (\widehat\E)^\vee \ra
\E xt^1(\G_2,\cO_{\PP'}) \ra 0;$$
$$0 \ra \E xt^m(\widehat\E,\cO_{\PP'}) \ra
\E xt^{m+1}(\G_2,\cO_{\PP'})\ra0,\ \forall m\geq 1.$$
We have now that $\pd \G_2=1$ if and only if
   $\E xt^m(\G_2,\cO_{\PP'})=0$, $\forall m>1$.
Thus $\E xt^m(\widehat\E,\cO_{\PP'})=0$ $\forall m>0$,
equivalently $\pd(\widehat\E)=0$ and $\widehat\E$ is locally free.

\qed

\bigskip
We are now going to describe the case where $d_1=3$, $d_2=d_3=2$,
and relate the above constructions (considering also
all the possible permutations of the spaces $W_1, W_2,W_3$)
to the geometry of cubic surfaces in $\PP^3$.
We consider now a tensor
\begin{equation}
     \widehat B\in V\otimes (\widehat U)^\vee \otimes \widehat{W};
\quad \dim V=4, \dim \widehat U=\dim \widehat W=3.
\end{equation}

Observe that we have 6 permutations of the three vector spaces,
inducing 6 distinct product projections. Moreover, we may vary the
twisting factor $t$.

\smallskip
We consider the exact order  given above of the three vector spaces.
On $\PP^3:=\PP(V^\vee)=\Proj(Sym(V))$,
for $t=1$ we get the sheaf $\widehat\G_1=({\G_1})_{\widehat B}$:
\begin{equation}
     0 \ra \widehat W ^\vee\otimes\cO_{\PP^3}(-1)
     \rxar{\widehat B_{\widehat W^\vee,\widehat U^\vee}}
     \widehat U^\vee\otimes\cO_{\PP^3} \ra \widehat\G_1 \ra 0\label{G1}
\end{equation}

For $t=2$ we get a vector bundle $\widehat \E:=\widehat\E_{\widehat B}$,
fitting in the two exact sequences
\begin{gather}
\label{Ehat}
0 \ra \wedge^2(\widehat W ^\vee) \otimes\cO_{\PP^3}(-2) \ra
\widehat W^\vee\otimes \widehat U^\vee\otimes \cO_{\PP^3}(-1) \ra
\widehat\E \ra 0\\
0 \ra \widehat\E \ra S^2 (\widehat U^\vee)\otimes \cO_{\PP^3}\ra
\widehat\G_2 \ra 0.
\end{gather}

For $t\geq 3$ it is pointless to proceed further, since indeed one
finds that $\widehat G_t$ is the $t$-th symmetric power of $\widehat G_1$.
In fact, by Corollary \ref{cor2}, $\widehat G_t$ has resolution:
\begin{equation}\label{Gt}
0\ra
\begin{matrix}
\wedge^{3}\widehat W^\vee \\
\otimes \\ S^{t-3}\widehat U \\
\otimes \\ \cO_{\PP^3}(-3)
\end{matrix}
\ra
\begin{matrix}
     \wedge^2\widehat W^\vee\\
\otimes \\ S^{t-2}\widehat U\\
\otimes \\ \cO_{\PP^3}(-2)
\end{matrix}
\ra
\begin{matrix}
     \widehat W^\vee\\
\otimes \\ S^{t-1}\widehat U\\
\otimes \\ \cO_{\PP^3}(-1)
\end{matrix}
\ra
\begin{matrix}
S^{t}\widehat U\\ \otimes \\ \cO_{\PP^3}
\end{matrix}
\ra \widehat\G_t \ra 0,
\end{equation}
which is the third symmetric power of (\ref{G1}).

\smallskip
We now consider the order
$V\otimes \widehat{W}\otimes (\widehat U)^\vee$
and repeat the same construction:
this is equivalent to consider the above construction
applied to the tensor
$$\sigma(\widehat B)\in V\otimes \widehat{W}\otimes (\widehat U)^\vee,$$
where $\sigma$ is the involution permuting $\widehat{W}$ with
$(\widehat U)^\vee$.

We obtain the sheaves $(\widehat\G_1)_{\sigma(\widehat B)}$,
$(\widehat\G_2)_{\sigma(\widehat B)}$, and
another vector bundle $\widehat\E_{\sigma(\widehat B)}$,
sitting in exact sequences equals to (\ref{G1})--(\ref{Gt})
with the roles of $\widehat{W}$ and $(\widehat U)^\vee$ exchanged:
\begin{gather}
     0 \ra \widehat U \otimes\cO_{\PP^3}(-1)
     \rxar{\widehat B_{\widehat U,\widehat W}}
     \widehat W\otimes\cO_{\PP^3} \ra (\widehat\G_1)_{\sigma(\widehat
     B)} \ra 0\\
0 \ra \wedge^2 \widehat U\otimes\cO_{\PP^3}(-2) \ra
\widehat U \otimes \widehat W  \otimes \cO_{\PP^3}(-1) \ra
\hat\E_{\sigma(\widehat B)} \ra 0\\
0 \ra \hat\E_{\sigma(\widehat B)} \ra S^2 \widehat W \otimes \cO_{\PP^3}\ra
\hat{\G}_2 \ra 0.
\end{gather}

As we shall see, these two sheaves,
supported on the same cubic surface
$G=\{\det(\widehat B_{\widehat W^\vee,\widehat U^\vee})=0\}
(=\{\det(\widehat B_{\widehat U,\widehat W})=0\})$,
correspond to two plane representations of $G$ as the
blow-down of a sixtuple of lines in a Double-6 configuration.
One passes from one plane representation to the other
by excanging the roles of $\widehat{W}$ and $(\widehat U)^\vee$,
i.e., by applying the trivial involution $\sigma$
to the tensor $\widehat B$.

\medskip
Considering the ordering
$\widehat B\in \widehat U^\vee\otimes V \otimes \widehat W$,
similar in spirit to the one
$\widehat B\in \widehat W \otimes V \otimes \widehat U^\vee$
we obtain a different geometric picture.
Recall that $\PP^3=\PP(V^\vee)$.
We are considering the projection
$p:\ \PP(\widehat U)\times\PP^3 \ra \PP(\widehat U)$ and
$\Gamma\subset \PP(\widehat U)\times\PP^3$
is the graph of the contraction
morphism from $G\subset \PP^3$ to $\PP(\widehat U)$.

For $t=1$, corollary (\ref{cor2}) provides the resolution
$$
0 \ra \widehat W^\vee\otimes\cO_{\PP(\widehat U)}(-1) \ra
V^\vee\otimes\cO_{\PP(\widehat U)} \ra p_*(\cO_G(1)) \ra 0,
$$
which is the Hilbert--Burch resolution of
$\I_\zeta(3)$, a twist of the ideal sheaf of a
length 6 0-dimensional subscheme $\zeta:=p(\Gamma)$ of $\PP(\widehat U)$.
Thus  $p_*(\cO_G(1))=\I_\zeta(3)$, that is
the linear forms on $G$ correspond to the cubics in $\PP(\widehat U)$
which are in the ideal sheaf $\I_\zeta$.

For $t=2$ we get a resolution of $p_*(\cO_G(2))$
{\small 
\begin{equation}
\wedge^2\widehat W^\vee\otimes\cO_{\PP(\widehat U)}(-2) \ra
\widehat W^\vee\otimes V\otimes\cO_{\PP(\widehat U)}(-1)
\ra S^2V\otimes\cO_{\PP(\widehat U)}\ra p_*(\cO_G(2)),
\end{equation}
}
and we find again the symmetric square of the previous resolution,
thus a resolution for ${\I_\zeta}^2(6)$.

Similarly for the cases where $t\geq2$.

\medskip
Quite interesting is instead the ordering
$\widehat B\in \widehat U^\vee\otimes \widehat W \otimes V^\vee$.
similar in spirit to the one
$\widehat B\in \widehat W\otimes \widehat U^\vee\otimes V^\vee$,
In this case
$p(\Gamma)=\{u\in \widehat U^\vee \mid \widehat B\neg u
\text{ has a kernel}\}$
and
$\Gamma\subset \PP(\widehat U) \times \PP(\widehat W^\vee)$
is the complete intersection of 4 hypersurfaces of bidegree $(1,1)$.
Let $H_1$ be the hyperplane class in  $\PP(\widehat U)$
and $H_2$ the one in $\PP(\widehat W^\vee)$:
since $(H_1+H_2)^4=6{H_1}^2{H_2}^2$,
we conclude that in general $\Gamma$ consists of 6 points
(but for our purposes it suffices
that $\Gamma$ is a complete intersection  0-dimensional subscheme
of length equal to 6).
Let again be $\zeta$ the length-6 subscheme $\zeta:=p(\Gamma)$.

For $t=2$ corollary (\ref{cor2}) is still applicable
and we get a non-classical resolution for
$\cO_\zeta=p_*\cO_\Gamma(0,2)$:
$$
0 \ra \wedge^2 V\otimes\cO_{\PP(\widehat U)}(-2)
\ra V\otimes W \otimes\cO_{\PP(\widehat U)}(-1)
\ra S^2 W \otimes \cO_{\PP(\widehat U)} \ra p_*\cO_\Gamma(0,2) \ra 0.
$$

For $t=1$, the complex
$$
V\otimes \cO_{\PP(\widehat U)}(-1) \ra
\widehat W\otimes\cO_{\PP(\widehat U)} \ra \cO_\zeta \ra 0
$$
is no longer necessarily exact, in the sense that corollary (\ref{cor2})
does not apply.
We shall now show  (cf. next corollary) that
we get
$$
0 \ra \cO_{\PP(\widehat U)}(-4) \ra V\otimes\cO_{\PP(\widehat U)}(-1)
\ra \widehat W\otimes\cO_{\PP(\widehat U)} \ra \cO_\zeta \ra 0.
$$
This exact sequence may also be obtained,
using
$\cO_\zeta\cong \E xt^2(\cO_\zeta,\cO_{\PP(\widehat U)}(-4))\cong
\E xt^2(\cO_\zeta,\cO_{\PP(\widehat U)})$,
as the dual of the Hilbert--Burch resolution of $\cO_\zeta$:
$$
0 \ra \widehat W^\vee \otimes\cO_{\PP(\widehat U)}(-4) \ra
V^\vee\otimes\cO_{\PP(\widehat U)}(-3) \ra \cO_{\PP(\widehat U)}
\ra \cO_\zeta \ra 0.
$$
The following corollary  spells out in detail  corollary \ref{cor1}.
\begin{cor}\label{cor3}
Suppose that $\Gamma$ is not empty, c.f. (\ref{Gammanotempty}).
Assume that $t>0$ but $t<d_3+1-d_2$.

Then $\G_t$ has a resolution given by an
Eagon--Northcott type complex:
{\small
\begin{equation}
\begin{split}
0\ra
\begin{matrix}
\wedge^{d_3+1}W_3^\vee \\
\otimes \\ S^{d_3-d_2-t}W_2^\vee \\
\otimes \\ \cO_{\PP'}(-d_3-1)
\end{matrix}
\ra \dots \ra
\begin{matrix}
\wedge^{t+d_2+2}W_3^\vee \\
\otimes \\ W_2^\vee \\
\otimes \\ \cO_{\PP'}(-t-d_2-2)
\end{matrix} \ra
\begin{matrix}
\wedge^{t+d_2+1}W_3^\vee \\
\otimes \\ \cO_{\PP'}(-t-d_2-1)
\end{matrix}
\ra \\
\ra
\begin{matrix}
\wedge^{t}W_3^\vee \\
\otimes \\ \cO_{\PP'}(-t)
\end{matrix}
\ra \dots \ra
\begin{matrix}
\wedge^1  W_3^\vee\\
\otimes \\ S^{t-1}W_2\\
\otimes \\ \cO_{\PP'}(-1)
\end{matrix}
\ra
\begin{matrix}
S^{t}W_2\\ \otimes \\ \cO_{\PP'}
\end{matrix}
\ra \G_t \ra 0.
\end{split}
\end{equation}
}
\end{cor}

\proof
Of course, we have
$$
H^j(\cO_{\PP^{d_2}}(a)))=0 \text{ except for }
\begin{cases}
     a\geq0 &\text{ if } j=0;\\
     a\leq -d_2-1 &\text{ if } j=d_2.
\end{cases}.
$$
Suppose that we have a free resolution on $\PP$ of $\cO_\Gamma(t)$
with terms
$$
0\ra L_r \ra \dots \ra L_{t+d_2+1} \ra \dots \ra L_t \ra
\dots \ra L_0 \ra \cO_\Gamma(t),$$
with $\deg L_j=t-j$.
By applying the functor
$ p_*(-)$ to an injective resolution,
we obtain a double complex whose vertical spectral term has
an $E_1$ term of the form
{\small
\begin{equation*}
\begin{matrix}
R^{d_2}p_*(L_r) &*\dots * &R^{d_2}p_*(L_{t+d_2+1}) &\dots &0 &\dots &0
\\
\dots &\dots &\dots  &\dots &\dots &\dots \\
\\
0 &\dots &\dots &\dots &R^0p_*(L_t) &*\dots* &R^0p_*(L_0)
\end{matrix},
\end{equation*}
}
whose non-zero terms are only the ones
indicated explicitly or with $*\dots*$.

Moreover $R^{d_2}p_*(L_{t+a})=
\wedge^{t+a}W_3^\vee\otimes
H^{d_2}(\PP(W_2),\cO(-a))\otimes \cO_{\PP'}(-t-a)$
and
$H^{d_2}(\PP(W_2),\cO(-a))\cong H^0(\PP(W_2),\cO(-d_2-1+a))$
by Serre's duality.
As proven in corollary \ref{cor1}, we obtain therefore
   complex, which is exact by lemma
(\ref{t>0}) since $t>0$, and is a resolution of  $\G_t$.
\qed

\section{Vector bundles $\E$ on $\PP^3$ with Chern polynomial
     $1+3t+6t^2+4t^3$:
general features, their construction, and cubic surfaces}


In this section we study some general features of vector
bundles $\E$ on $\PP^3$ with Chern polynomial $1+3t+6t^2+4t^3$.
Recall lemma \ref{E=beta} of section 1:
under the open condition of having minimal
cohomology, these bundles
have quite a simple copresentation
in terms of their intermediate cohomology modules.
Indeed, we have seen that $H^2_*(\PP^3,\E)=0$ and
that $H^1_*(\PP^3,\E)$ has only two nonzero graded pieces,
namely the vector spaces $U=H^1(\PP^3,\E(-2))$ and
$W=H^1(\PP^3,\E(-1))$.
Recall moreover that $V:=H^0(\Pn 3,\cal O(1))$
is the space of linear forms on $\Pn3$.

We will see that there are three ways
to construct such bundles:
\begin{enumerate}
\item as syzygy bundles
starting from a tensor
\begin{equation}  \label{Btensor}
B\in U ^{\vee}\otimes W\otimes V^{\vee},
\end{equation}
which will be our natural choice to parametrize $\E$
(we shall call this the {\bf kernel construction});
\item   as  extensions,
starting from another tensor
$ B' \in {W'}^\vee \otimes {U'}^\vee \otimes  V $
(we shall call this the {\bf direct construction})
\item as a {\bf direct image sheaf},
starting from a third tensor 
$\widehat B\in V\otimes {\widehat U}^\vee \otimes {\widehat W}$
and using the construction described in section 2.
\end{enumerate}
The relation occurring between $B$ and $ B'$ will lead to the definition
of the cross-product involution,
while the relation occurring between $B$ and $\widehat B$
will be investigated in the section 4
after we introduce the cross-product
involution.
\bigskip

We now explain the first construction.
Suppose we have  a bundle $\E$ as in (\ref{E=ker}).
Applying $\Hom(-,\cal O)$ to the Euler sequence (\ref{euler})
and tensoring by $\hol(1)$ yields,
since $\Hom(\hol(2),\hol(1))=0$, $Ext^1(\hol(2),\hol(1))=0$,
\begin{equation}\label{B=beta}
\Hom(\Omega^1(2),\cal
O(1))\cong \Hom( V \otimes \cal O(1) , \cal O(1)).
\end{equation}

Thus the map $\beta$ factors through a map
$B: U\otimes(V \otimes\hol(1)) \rar W\otimes\hol(1)$ and
the sheaf map $B$ is surjective.
This surjectivity is obviously equivalent to
$H^0(B(-1)):\ U\otimes V \ra W$ being surjective.

In the sequel we shall often identify the sheaf map $B$ with the
corresponding tensor
$H^0(B(-1)) \in U ^{\vee}\otimes W\otimes V^{\vee}$.

Let  $\epsilon$ be the tensor product of the identity map of $U$
with the evaluation map
$V\otimes \cal O\rar \cal O(1)$.
Then one sees easily that
$\cal E =\ker \beta=\ker B\cap\ker \epsilon$
and that the short exact sequence (\ref{E=ker}) can be replaced by
\begin{equation}\label{E=ker2}
          0 \rar \cal E \rar U\otimes V \otimes\cal O(1)
\rxar{B\oplus\epsilon}
(W\otimes\cal O(1))\oplus (U\otimes\cal O(2)) \rar 0,
\end{equation}

\begin{rem}
The cohomology  exact sequence associated to the following
twist of (\ref{E=ker2}), namely:
$$ 0 \rar \cal E(-2) \rar U\otimes V\otimes \hol(-1)
\rxar{(B\oplus\epsilon)(-2)}
(W\otimes \hol(-1)) \oplus (U\otimes \hol) \rar 0$$
yields a canonical isomorphism $ U\cong H^1 (\cal E(-2))$.

Since there is a canonical isomorphism
$$H^0(\epsilon(-1)): U\otimes V \ra U\otimes H^0(\hol(1)),$$
the projection of
$W\oplus(U\otimes V) \ra W$ induces an isomorphism of
$H^1(\cal E(-1))= Coker H^0((B \oplus \epsilon)(-1))$ with
$W$, such that the map
$B:\ U\otimes V \rar W$
corresponds to the multiplication map of the
cohomology module $ H^1_{*} (\E)$.
\end{rem}

\begin{defn}
The {\bf  kernel construction} of the bundle $\E$ is
as follows.
Consider a $3 \times 3 \times 4$ tensor
\begin{equation}
B\in U^\vee\otimes W \otimes V^\vee.
\end{equation}

Such $B$ induces a linear map $B:  V \otimes U \ra W $
and a homomorphism
$B:\  V \otimes U \otimes \hol \ra W \otimes \hol$
of vector bundles on $\PP^3:=\PP(V^\vee)$,
which induces a homomorphism
$\beta=B\oplus\epsilon: U \otimes \Omega^1 (2) \ra  W \otimes \hol(1)$
as described above.

If $\beta$ is surjective,
$\E : = ker (\beta)$ as in (\ref{E=ker})
is a vector bundle.
\end{defn}

Moreover, lemma \ref{E=beta} shows that such an $\E$ is a vector
bundle with total Chern class  $c(\E) (t) = 1 + 3 t + 6 t^2 + 4 t^3$
and, if $H^0(\beta)$ is surjective, with minimal cohomology
(i.e., the conditions (1)--(5) of lemma \ref{E=beta} are satisfied),
and, moreover, $U=H^1(\PP^3,\E(-2))$, $W=H^1(\PP^3,\E(-1))$ and
the multiplication tensor is exactly $B$.

\bigskip
We now proceed by illustrating the direct construction, 
our second construction.
Recall again that, by lemma \ref{E=beta}, 
a tensor  $B$ such that $H^0(B\oplus \epsilon)$ is surjective 
gives an $\E$ with $H^i(\E) = 0 \ \forall i \geq 1$, and therefore
$h^0(\E)=6$.
It seems therefore  natural to introduce the so-called

\medskip\noindent
{\bf SECOND ASSUMPTION: }\\
1) $ \i : 6\cal O \rar \cal E $ {\bf is injective},
hence  we get an exact sequence:
\begin{equation}\label{seq1a}
               0\rar 6\cal O \rar \cal E \rar \tau \rar 0,
             \end{equation}
{\bf 2)
the torsion sheaf $\tau$  is $\cal O_G$-invertible,} where $G$ is the
divisor of $\Lambda^6(\i)$.
\medskip

If for a vector bundle given as in (\ref{E=ker})
the second assumption is satisfied, then the divisor
$G$ is a cubic surface, and $\E$ may be reconstructed
as an extension  of $6\cO$ by $\tau$.

\smallskip
We first analyse the geometry and cohomology of the $\hol_G$-invertible
sheaf $\tau$.
Let $\tau=\cO_G(D)$ and let $H$ denote
the hyperplane divisor in $\PP^3$.
We reproduce here remark 4.3 of \cite{CaTo}

\begin{rem}\label{tau}
Notation being as 1) above (even without assuming $\tau$ to be
$\hol_G$-invertible) set $\tau'=\cal Ext^1(\tau,\cal O)$:
then the dual of the previous exact sequence
(\ref{seq1a}) gives
             \begin{equation}\label{seq1b}
0\rar \cal E^\lor \rar 6\cal O \rar \tau' \rar 0.
             \end{equation}
and we have:
\begin{enumerate}
\item By (\ref{seq1a}) clearly  $H^0(\tau)=H^1(\tau)=H^2(\tau)=0$.

\item By (\ref{seq1b}) and since $h^i(\cal E^\lor)\cong h^{3-i}(\cal E(-4))$ we
get
            $h^0(\tau')=6$, $H^1(\tau')=H^2(\tau')=0$.
\item Since by definition $\tau'=\cal Ext^1(\tau, \cal O)$, applying
the functor
$\cal Hom(\tau,-)$ to the exact sequence
$0\rar \cal O\rar \cal O(3) \rar \cal O_G(3)\rar 0$ we get
$\tau'=\cal Hom (\tau,\cal O_G(3))$.
Therefore, if $\tau=\cal O_G(D)$, then $\tau'=\cal O_G(3H-D)$.

\end{enumerate}
\end{rem}
Since  $h^i(D)=0\  \forall i$,
$h^0(3H-D)=6$, $h^i(3H-D)=0$
for $i=1,2$,  by Riemann Roch follows that
$D^2+DH=-2$ and $10=36-7DH+D^2$. Therefore $HD=3$, $D^2=-5$.

Setting $\Delta:=D+H$, it turns out that $\Delta
H=6,\Delta^2+\Delta K_G=-2$, i.e.,
$$\Delta H=6,\Delta^2=4.$$

\begin{lem}\label{D}
Assume that $G$ is a smooth cubic surface: then
there exists a realization of  $G$ as the image of the plane
under the system
$| 3L - \sum_1^6 E_i|$ of plane cubics through six points, such that either
$\Delta \equiv 2L$, i.e., $\Delta$
corresponds to the conics in the plane, or (up to permutations of the six
points)
$\Delta \equiv 3L - 2 E_1 - E_2.$
\end{lem}

\Proof
Observe preliminarly that if $ |H| = | 3L - \sum_1^6 E_i|$
is such a planar realization of a cubic surface, then another one
is obtained via a standard Cremona transformation
centered at three of the points $P_i$ corresponding to the
$(-1)$-curves $E_i$.

In fact, if $ L' : = 2L - E_1 - E_2 - E_3 $, then
$$ H =  3L'  - (L - E_1 - E_2) -  (L - E_1 - E_3) - (L - E_2 - E_3)
- E_4 - E_5 - E_6. $$

We have $0= H^2(D)= H^0(-D-H)$ and
     a fortiori $H^2(\Delta)=H^0(-D-2H)=0$.
It follows that
$|\Delta|$ has $h^0(\Delta)\geq 6$, $\Delta^2=4$, and the arithmetic
genus $ p_a(\Delta)=0$.

Since $ 0 = H^1(D)= H^1( K_G-D) = H^1(-H -D) = H^1(- \Delta)$,
     it follows that $\Delta$ is connected.

Hence we have a representation
$\Delta \equiv n L - \sum_1^6  a_i E_i$, where the $a_i$'s are non negative
and we assume $a_1\geq a_2\geq \dots\geq a_6$.

We have: $ \Delta^2 = 4 = n^2 - \sum_1^6  a_i^2  $, $\Delta
\cdot H=6 = 3 n - \sum_1^6  a_i$, i.e.
          \begin{equation}\label{nn2}
n^2 = \sum_1^6  a_i^2 + 4 , \  3 n = \sum_1^6  a_i + 6 .
          \end{equation}

We want to show that, after a suitable
sequence of standard Cremona transformations,
$\Delta\equiv 2L$  or $\Delta\equiv 3L-2E_1-E_2$.
By (\ref{nn2}), we have $n\geq 2$
and for $n=2,3$ $\Delta$ is as claimed.
Hence the claim is that there exists a
sequence of standard Cremona transformations
which makes  $|\Delta|$
have degree $n\leq 3$.

By applying $|2L-E_1-E_2-E_3|$ we get a new system $\Delta'$
with degree $n'=2n-a_1-a_2-a_3$.

By our ordering choice for the $a_i$'s, we have
$$a_1+a_2+a_3 \geq (\sum_{i=1}^6 a_i)/2=3n/2 -3,$$
with strict inequality unless all $a_i$'s are equal.
We study this latter case  first:

\bigskip
{\bf Sublemma.} In the previous setting, $a_1=a_2=\ldots=a_6$ if and only
if $n=2$ and $a_1=a_2=\ldots=a_6=0$ or $n=10$ and $a_1=a_2=\ldots=a_6=4$.

\medskip
\proof
The statement follows immediately by defining $a:=a_1=a_2=\ldots=a_6$ and
using both conditions of (\ref{nn2}): $n=2a+2$, $n^2=6a^2+4$,
which imply $ 8 a =2 a^2$.
\qed

\bigskip
The previous inequality gives:
$$n'\leq\frac{n}2+3\leq n \text{\qquad for $n\geq 6$},$$
and $n'<n$ for $n\geq 6$ unless $n'=n=6$ and $a_1=a_2=\ldots=a_6$,
which has no solution by the above sublemma.
We conclude that after  suitable Cremona transformations $n\leq 5$.

If $n=5$, then $n'\leq 5/2+3$, i.e., $n'\leq 5$.
Moreover, if also $n'=5$, then $a_1+a_2+a_3=5$ and using again
\ref{nn2} we obtain $a_4+a_5+a_6=4$. But then $a_1=a_6+1$
and we easily get
a contradiction since then $a_2= a_3=a_4=a_5 $
and they either equal $a_1$ or $a_6$.
Hence, after a suitable Cremona transformation,
we can always reduce to the case $n\leq 4$.

Let now $n=4$. Using (\ref{nn2}) we get
$\sum_{i=1}^6 a_i=6$ and $\sum_{i=1}^6 a_i(a_i-1)=6$.
We have the following two possibilities:
$a_1=3, a_2=a_3=a_4=1, a_5=a_6 =0$ or $a_1=a_2=a_3=2,
a_4=a_5= a_6 =0$.
In both cases we have that $a_1+a_2+a_3\geq 5$, and therefore $n'\leq 3$.

%

\QED

\begin{rem}
The complete linear system $\Delta$ has as image in $\PP^5$ either
the Veronese embedding of $\PP^2$, or the embedding
of $\PP^1 \times \PP^1$ through
$ H^0 ( \hol_{ \PP^1 \times \PP^1 } (1,2))$.
In both case we have a surface of minimal degree (=4).
\end{rem}

\bigskip
Thus we have concluded that either $D=2L-H=-L+\sum E_i$, or
$D=3L- 2 E_1 - E_2 - H=- E_1 + \sum_3^6 E_i$.
The latter case does not occur, because an extension
of   $6\cO$ by such a $\tau$ will not have the required cohomology
table.
The former case is instead  possible.
We refer to \cite{CaTo}, Lemma 4.12 and Lemma 4.13 for the proof
of these facts.

\medskip
We are now able to explain the second way
to construct vector bundles $\E$
such that $ H^i(\E) = 0 \ \forall i \geq 1$ and $h^0(\E) = 6$,
as required in Lemma (\ref{E=beta}):
we construct them as  extensions of $6\cO$ by $\tau$,
where
$\tau$ is the sheaf corresponding to $-L+\sum_1^6E_i$
on a smooth cubic surface $G$ (cf. \cite[Lemma 4.12]{CaTo}).

Before we give this second construction,
we study these extensions.
Setting $\tau' : = \E xt^1(\tau,\cal O)$ and recalling
remark (\ref{tau}), we see that
such  extensions are parametrized by
$Ext^1(\tau,6 \cal O)=H^0(6\cal Ext^1(\tau,\cal O))\cong \CZ^{36}$.

\begin{lem}\label{unicity}\cite[Lemma 4.10]{CaTo}
Assume that $h^0 (\E^\vee ) = 0$ and that $\E$ is an
extension  as in (\ref{seq1a}): 
then the extension class in
$Ext^1(\tau,6 \cal O)=H^0(6\cal Ext^1(\tau,\cal O))\cong \C^{6}
\otimes \C^6$ is a rank 6 tensor
(we shall refer to this statement by saying that
{\bf the extension does not partially split}).

In particular,
{\bf $\E$ is then uniquely determined up to  isomorphism}.
\end{lem}
\proof

The extensions which yield vector bundles
form an open set.

We canonically view the space of these extension classes
as $ {\rm Hom}(H^0(\tau'),H^0( 6 \cal O)) = {\rm Hom} (H^0(\tau'),\C^6)$,
through the coboundary map of the corresponding exact sequence.
We have then an action of $GL(6, \C)$ as a group of automorphisms of
$6 \cal O$, which induces an action on $ {\rm Hom}(H^0(\tau'),H^0( 6 \cal
O)) = {\rm Hom}(H^0(\tau'),\C^ 6 ) $ which corresponds to
the composition of the corresponding linear maps.

The extensions which yield vector bundles
form an open set, which contains an open dense orbit,
on which this action is free, namely, the tensors of rank = 6.

If the rank of the tensor corresponding to an extension is $ = r < 6$,
it follows that the extension is obtained from an extension
$ 0 \ra r \hol \ra \E'' \ra \tau \ra 0 $ taking then a direct sum
with $ ( 6 - r) \hol$: but then $ ( 6 - r) \hol$  is a direct summand
of $\E^\vee$,  contradicting $h^0 (\E^\vee ) = 0$.

\qed

\begin{cor}\label{vb}
$\E$ as above (\ref{unicity}) is a vector bundle if
$H^0(\tau')$ has no base points.
\end{cor}

\proof
Our hypothesis shows that in each point of $G$ the local extension class is
non zero, hence it yields a locally free sheaf.
\qed

\bigskip

The second assumption yields a cubic surface $G \subset \Proj(V)$
and an invertible sheaf $\tau$ on $G$. If $G$ is smooth,
the invertible sheaf
$ \tau (1)=2L$ yields then a birational morphism onto
a Veronese surface, whence represents $G$ as the blow up of a projective
plane $\PP^2$  in a subscheme $\zeta$ consisting of six points (distinct
if the cubic $G$ is smooth), and
as the image of
$\PP^2$ through the linear system of cubic curves  passing through
$\zeta$.
The Hilbert-Burch theorem allows us to make an explicit construction
which goes in the opposite direction.

\begin{rem}
Let $U',W'$ be 3-dimensional vector spaces and set $\PP^2 :=
\PP ({U'})$.
Consider a $3 \times 3 \times 4$ tensor
\begin{equation}\label{hatb1}
     B'\in {W'}^\vee \otimes {U'}^\vee \otimes V
\end{equation}
and assume that the induced sheaf homomorphism
$W'\otimes \cal O_{\PP^2}(-1) \ra V\otimes \cal O_{\PP^2}$,
which we call again $B'$,
yields an exact sequence
\begin{equation}
0 \rar W'\otimes \cal O_{\PP^2}(-1) \rxar{B'}
V\otimes \cal O_{\PP^2}
\rxar{ \Lambda^3({B'}) } \cal O_{\PP^2}(3) \rar \cal O_\zeta(3)
\rar 0
\end{equation}
which is the Hilbert Burch resolution of a codimension $2$ subscheme
$\zeta$
of length $6$.

We obtain a canonical isomorphism $V\cong H^0( I_\zeta(3))$
and we let  $G\subset\Proj(V)$
be the image of $\PP^2$ via the rational map $\psi$
associated to $V$.
Under the above assumption on  $B'$,
if moreover $\zeta$ is a local complete intersection,
$G$ is a normal cubic surface and, if we set
$\cal G:=(\psi_* (\cal O(1))$,
   there is an  exact sequence on $\Proj(V)$:
\begin{equation}\label{Geq1}
0 \rar W'\otimes \cal O (-1)  \rxar{B'}
{U'}^\vee\otimes \cal O \rar \cal G \rar 0.
\end{equation}
\end{rem}

Under the more general assumption that the sheaf map $B'$
in (\ref{Geq1})
never drops rank by 2, $\cal G$ is an invertible
sheaf on the cubic surface $G$, and there is a Cartier divisor $L$
on $G$ such that $\cal G=\cal O_G(L)$ 
(and $h^0(\cal O_G(L))=3$).

\begin{defn}\label{dir_cos}
We define now the {\bf  direct construction } of the bundle $\E$
relying on our results above.

Consider a $3 \times 3 \times 4$ tensor
\begin{equation}\label{hatb}
B'\in {W'}^\vee \otimes {U'}^\vee \otimes V
\end{equation}
such that
the sheaf $\cal G$ defined by the exact sequence (\ref{Geq1})
is invertible on the cubic surface $G$ (i.e., at each point
$ y \in \PP^3$ $ rank ( \cal G \otimes \C_y ) \leq 1$).

Define  $ \tau : = \cal G^{\otimes2}(-1)$ and let $\E$ be
a vector bundle which is an extension of $ 6 \hol$ by $\tau$
         as in (\ref{seq1a})
(here and elsewhere, $\hol : = \hol_{\Proj(V)}$).
\end{defn}

We then have the following results:

\begin{prop}\label{un2}\cite[Proposition 4.17]{CaTo}
         $\E$ as above is unique up to isomorphism in the following cases:
\begin{enumerate}
\item
         if $G$ is a smooth cubic surface.

\item
          if $G$ is  the reduced union of a plane $T$ and
a smooth quadric $Q$ intersecting transversally.
\end{enumerate}
\end{prop}

\begin{rem}\label{doubleline}
The case where $G$ is a a linear projection of the cubic scroll $Y$
   (birational embedding of $\PP^2$ in $\PP^4$ by the system $ | 2 L - E|$)
yields two sheaves $ \cal G$ which are not invertible.

As it is well known, every point in $\PP^4$ lies in one of the planes
spanned by the conics of the system $L$.
If we project from a point in $\PP^4\setminus Y$, this conic maps two to
one to the double line of the cubic $G$.

Such a plane is said to be {\bf special} if the conic splits
into two lines $ E + F, F \equiv L-E$.

In the non special case,  we may assume without loss of generality
that the conic corresponds to the line $ z = 0$ in the plane,
that the blown up point is the point $ x=y=0$, and that
the linear system mapping to $G$ is generated by
$ ( zx : =  x_0,   zy :  = x_1, x^2 : = x_2, y^2 : = x_3) $.
In this case one sees that the matrix $B'$ is
$$B'=\begin{pmatrix}
x_0 &0 &x_1\\ 0 &x_1 &-x_0\\ -x_2 &-x_3 &0
\end{pmatrix},$$
hence the rank of $B'$ drops by 2 on the
line $ x_0= x_1 = 0$
($G$ is then the cubic of equation
$-x_0^2x_3+x_2x_1^2=0$).

In the special case, we may again assume that
the blown up point is the point $x=y=0$,
we assume that the line $F$ is the proper transform of $x=0$,
and that the linear system mapping to $G$ is generated by
$ ( y^2+zx : =  x_0,   x^2 :  = x_1, yz : = x_2, xy : = x_3) $
( in the projective embedding given by $ (zx, yz, x^2, xy, y^2)$
it corresponds to projection from the point $(1,0,0,0,-1) \in \PP^4
\setminus Y$).

In this case one sees that the matrix $B'$ is
$$B'=\begin{pmatrix}
x_3 &-x_2 &-x_0\\ -x_1 &0 & x_3\\ 0 &x_3 &-x_1
\end{pmatrix},$$
hence the rank of $B'$ drops by 2 on the
line $ x_3= x_1 = 0$ ($G$ is then the cubic of equation
$-x_3^3+x_1^2x_2+x_1x_2x_3=0$).
\end{rem}

\section{The cross-product-involution
and Schur's quadric}

In the previous section we have seen that to
a vector bundle as in (\ref{E=ker}) satisfying the second
assumption one can associate two tri-tensors:
the tri-tensor $B\in U^\vee\otimes W\otimes V^\vee$
and the tri-tensor
$B'\in {W'}^\vee \otimes {U'}^\vee \otimes V$.
The first corresponds to the unique nonzero multiplication matrix
of the intermediate cohomology module $H^1_*(\E)$, the
second, according to the direct construction,
defines on a cubic surface the invertible sheaf $\cal G$
such that $\E$ is an extension of $\tau=\cal G^{\otimes 2}(-1)$
and $6\cO$.

What is the relation between them?
In this section we will show that there is indeed a strict
relation between such tri-tensors:
a birational involution,
which the authors call {cross-product-involution}.

In \cite{CaTo} the authors, after having discovered these
two tensors, relate them by constructing a
not necessarily minimal resolution of a bundle $E$ constructed
by means of the tri-tensors $B$ and $B'$.

Indeed, given $B$, Beilinson's complex for $\E$
yields a short exact sequence
$$0 \rar U\otimes \Omega^2(2) \rar W\otimes
\Omega^1(1)\oplus 6\cal O \rar \cal E \rar 0,$$
where $U=H^1(\cal E(-2))$ and $W=H^1(\cal E(-1))$). We get:
{\small
\begin{equation}
          \label{eq:resE1}
0\rar
\begin{matrix}
          U\otimes \hol(-2)\\
          \oplus \\ 
          W\otimes \hol(-3)
\end{matrix}
\rar
\begin{matrix}
          U\otimes V\otimes \hol(-1)\\
          \oplus\\
          W\otimes \Lambda^3V\otimes \hol(-2)
\end{matrix} 
\rar
\begin{matrix}
          6\hol\\
          \oplus\\
          W\otimes\Lambda^2V\otimes\hol(-1)
\end{matrix}
\rar \E \rar 0.
\end{equation}
}

On the other side, $B'$ gives a resolution of $\cal G$,
from which it is possible to compute a resolution of
$\tau=\cal G^{\otimes 2}(-1)$.
From this one, by applying the mapping
cone, it is possible again to get a resolution of $\cal E$:
{\small
\begin{equation}\label{eq:resE}
0\rar \Lambda^2 W' \otimes\cal O(-3) \ra
{U'}^\vee \otimes W'  \otimes \cal O(-2) \ra
\begin{matrix}
  6\cal O \\
  \oplus \\
  S^2 {U'}^\vee  \otimes \cal O(-1)
\end{matrix}
\rar \cal E \rar 0.
\end{equation}
}

Comparing the two resolutions, the authors obtained the
following identifications:
\begin{enumerate}
\item 
$W'\cong \Lambda^2 W$,
$U\cong \ker[{\neg B'}:\ \Lambda^2 W'\otimes V^\vee \rar
{U'}^\vee\otimes W']$,
\item 
${U'}^\vee\otimes W' \cong ( W\otimes \Lambda^3 V)/U$,
$S^2{U'}^\vee\cong (W\otimes \Lambda^2V)/(U\otimes V)$;
\end{enumerate}
where ${\neg B'}$ is the contraction given by the composition
of the natural inclusion from
$(\Lambda^2 W') \otimes V^\vee$ to
$(W'\otimes W')  \otimes V^\vee$
with the map $B'\otimes id_{W'}(-1)$.

\smallskip
Based on the above considerations we state the following
correspondence for a pair of tri-tensors as above.

\begin{defprop}\label{cross}
\ \\
The {\bf Cross-Product Involution on Tensors of type 3 x 3 x 4}
is given as follows.

Consider a 4-tuple $(U,W,V^\vee,B)$,
where:
\begin{enumerate}
\item $U$, $W$ are 3-dimensional vector spaces and
$V$ is a 4-dimensional vector space;
\item each vector space is equipped with a given orientation,
identifying respectively
$\wedge^3 U$,  $\wedge^3 W$, and $\wedge^4(V^\vee)$ with $\CZ$;
\item $B$ is a tensor
$$B\in U^\vee \otimes  W\otimes V^\vee \cong \Hom(U,W)\otimes
V^\vee.$$
\end{enumerate}

We remark that the ordering of the three vector spaces yields
a sheaf homomorphism which is canonically associated
to the tensor $B$, namely
  $U\otimes \cO_{\Proj(V^\vee)}(-1) \rxar{[B]}W\otimes
\cO_{\Proj(V^\vee)}$.

\smallskip

The {\bf trivial involution} associates to the 4-tuple
$(U,W,V^\vee,B)$ the 4-tuple $(W^\vee,U^\vee,V^\vee,B)$.

The {\bf reversing construction} associates to the 4-tuple
$(U,W,V^\vee,B)$ the 4-tuple
$(W',{U'}^\vee,V,B')$,
where:
\begin{enumerate}
\item $W':=\Lambda^2(W)$. Since $W$ is equipped with
an orientation,
the duality
$W\otimes \Lambda^2 W \rar \C $ induces a canonical
isomorphism of
$W'$ with $W^\vee$.

\noindent
$U':=\ker[\neg B:\ \Lambda^2 (W^\vee)\otimes V \rar
U^\vee\otimes W^\vee ]$,
where ${\neg B}$ is the contraction with the tensor $B$;
in particular, $U'$ is canonically isomorphic to a
subspace of $W \otimes V$.

\item the three vector spaces $W'$, $U'$, and $V$ are equipped
with the orientations induced from the orientations of
$U$, $W$, $V^\vee$, under the 'main assumption' that the
contraction map  $\neg B$ be surjective.

\item the tensor
$$B'\in {W'}^\vee \otimes {U'}^\vee \otimes
V=\Hom(W',{U'}^\vee)\otimes V,$$ which corresponds
to the natural inclusion $ U' \ra W \otimes V$, in view of
the isomorphism $ W \cong ({W'})^\vee $.
\end{enumerate}

The dimension of $U'$ is equal to 3 if we make  the

{\bf MAIN ASSUMPTION: }
The contraction $\neg B$ is surjective.

The {\bf cross-product involution} is the involution,
defined for the 4-tuples  $(U,W,V^\vee,B)$, where $B$ is assumed to belong
to the open set of tensors
satisfying the main assumption, which is given by the composition of
the reversing construction with the trivial involution.

The cross-product involution associates thus
to the 4-tuple $(U, W, V^\vee, B)$ the 4-tuple
$(U', W = {W'}^\vee ,V ,B')$.

\end{defprop}

\begin{proof}
We only need to show that the cross-product involution
is an involution, i.e., that applying it twice, we obtain the identity.

We present here a different proof from the one given in \cite{CaTo},
and based on the following

{\bf Fact:} if the main assumption holds,
then there is an element $Q \in S^2 (V)$,
called Schur's quadric, such that if we denote by $ q : V^\vee \ra V$
the corresponding linear map, then $q$ is an isomorphism and
  $ id_W \otimes q$ carries $U \subset W \otimes V^\vee$
to $U' \subset W \otimes V$.

Indeed, by the construction of the Schur quadric,  it follows that the inverse
$q^{-1}$ of the linear map
$q$ is obtained from the Schur quadric $Q^\vee \in S^2 (V^\vee)$
associated to $B'$, and therefore $ U'' = id_W \otimes q^{-1} (U') = U$.

\end{proof}

\begin{rem}
The two tensors considered in remark \ref{doubleline},
whose respective determinants yield the two non normal irreducible cubics
(which are not projectively
equivalent) satisfy
the main assumption. But the cross-product involution constructs out
of them two tensors which do no longer satisfy the main assumption,
and which are projectively equivalent:

$$\begin{pmatrix}
x_2^\vee &x_3^\vee &0\\
0 &x_2^\vee &x_3^\vee\\
0&0&0
\end{pmatrix} \qquad
\begin{pmatrix}
0&0&0\\
x_0^\vee &x_2^\vee &0\\
x_2^\vee &0 &x_3^\vee
\end{pmatrix}.$$
\end{rem}
\bigskip \ \bigskip\

Let us now explain how the Schur quadric is obtained.

Let $U,W,V$ be complex vector spaces of respective dimensions 3,3,4
  and let $B\in U^\vee\otimes W\otimes V^\vee$
be a tensor of type $(3,3,4)$, as  in (\ref{Btensor}).
Following the notation of \cite{dk},
to $B$ are associated then 3 maps:
\begin{equation}\label{g}
g_{V^\vee}:\ V \ra U^\vee\otimes W,
\end{equation}
(we think of it as a $3\times 3$ matrix  of linear forms on $V$),
and similarly we view
$$g_{U^\vee}:\ U \ra W\otimes V^\vee,$$
$$g_W:\ W^\vee \ra U^\vee\otimes V^\vee,$$
as $3\times 4$ matrices of linear forms (respectively on $U$ and
$W^\vee$).

For a general $B$, the determinant of the $3\times 3$ matrix $g_{V^\vee}$
of linear forms on $V$ gives a smooth cubic surface  $G^*$
in the dual projective space ${\PP^3}^\vee = \Proj (V^\vee)$,
together with two different realizations of $G^*$
as a blow up of a projective plane $\Proj (U^\vee) $
(respectively $\Proj (W)$) in a set of six points $Z$.
These are the points where the $3\times 4$ Hilbert--Burch
matrix of linear forms on $U$ (respectively on $W^\vee$)
drops rank by 1,
and the rational map to ${\PP^3}^\vee$  is given
by the system of cubics through the 6 points,
a system which is generated by the determinants of the
four minors of order 3
of the Hilbert--Burch matrix,
the matrix $g_{U^\vee}$ (resp. $g_W$).
One passes from one realization to the other one
simply by applying the trivial involution to the tensor $B$, i.e.,
replacing $g_{U^\vee}$ with $g_W$.

Also the 12 lines of the double--six configuration can be
obtained  from the original tensor $B$,
as the union of the 6 lines $A_z=Ker(g_{U^\vee})$
with the 6 lines ${A'}_z=Ker(g_W)$
for $z\in Z$, cf. \cite[$\S$ 0]{dk}.
According to this notation,
Dolgachev and Kapranov give the following modern formulation
of Schur's classical theorem in \cite{Schur}:

\begin{thm}\cite[Theorem 0.5]{dk}
Given a smooth cubic there exists a
symmetric bilinear form $Q(x,y)$ on $V$,
unique up to a scalar factor, which satisfies the following
property:
$Q(x,y)=0$ whenever $x\in A_z$ and $y\in {A'}_z$ for some
$z\in Z$
(i.e., the corresponding lines of the double--six are
orthogonal with respect to $Q$). $Q$ is nondegenerate.
\end{thm}

The bilinear form $Q\in S^2(V)$ is called the  Schur quadric,
and it is obtained as follows.

Given a tri-tensor $B\in U^\vee\otimes W\otimes V^\vee$,
consider the second symmetric power of the linear map $ g = g_{V^\vee}$,
$$S^2(g) \colon \ra S^2(U^\vee\otimes W)$$
and compose it with the projection of
$$S^2(U^\vee\otimes W)=
\left(\Lambda^2 U^\vee\otimes \Lambda^2 W \right)
\oplus \left( S^2U^\vee \otimes S^2 W \right)$$
onto the first factor.

Since $\dim S^2 V=10$, $\dim \left(\Lambda^2 U^\vee\otimes \Lambda^2
W \right)=9$,  the kernel is 1-dimensional for a general tensor 
(cf. \cite[$\S0$ and Thm 0.5]{dk}).
\medskip

Recall once more that the cross-product involution associates to a 
general tensor $B\in  U^\vee \otimes W\otimes  V^\vee $ 
another tensor $B' \in {U'}^\vee \otimes W \otimes V$,
where $U'$ is defined as the kernel of the map
$\Lambda^2 (W^\vee)\otimes V \ra  U^\vee \otimes W^\vee $
induced by contraction with $B$.
\medskip

Associate to the  Schur quadric $Q \in S^2 V$ a linear map 
$q: V^\vee \ra V$. The map $q$ then relates $B$ and $B'$ as follows.

\begin{prop}
Let $B\in  U^\vee \otimes W\otimes  V^\vee $ be a tri-tensor
such that the associated cubic surface $G^* \subset {\PP^3}^\vee$
is smooth (in particular, $B$ and $B'$ lie in the open set of the tri-tensors
where the cross-product involution is defined).

Then the composition
of $g_{U^\vee} :\ U \ra W \otimes V^\vee$ with
\begin{equation}
     id_{W}\otimes q:\ W\otimes V^\vee \ra W\otimes V
\end{equation}
maps $U$ to $U'$, where $U'$ is the vector space associated
to $U$ via the cross-product involution.

In particular, the tensor $B'$, corresponding to the inclusion
  $ U'  \ra W \otimes V$, is determined in this way by the tensor $B$
and by the Schur quadric $Q$.
\end{prop}

\proof
According to the definition of the cross-product involution,
we can identify $\Lambda ^2 W^\vee$ with $W$,
$W^\vee$ with $\Lambda ^2 W$, and moreover

\begin{equation}\label{idcross}
\begin{aligned}
    U'=Ker(W\otimes V \ra U^\vee \otimes W^\vee),\\
U =Ker(W\otimes V^\vee \ra {U'}^\vee\otimes W^\vee), \\
\end{aligned}
\end{equation}
and both spaces have dimension equal to 3.

Therefore, since $q$ is invertible, in order to show that
$(id_{W}\otimes q )(g_{U^\vee}(U))=U'$, it suffices to show
that $(id_{W}\otimes q )(g_{U^\vee}(U))$ is contained in $U'$,
i.e., this space maps to zero in
$U^\vee \otimes W^\vee$.

Recall that the first map in (\ref{idcross}) is the composition
$$
W\otimes V \rxar{\otimes B}
(W\otimes W)\otimes (V\otimes V^\vee)\otimes U^\vee
\ra U^\vee \otimes W^\vee,
$$

where the second map is naturally obtained by the projection
$p_{W^\vee}:\ W\otimes W \ra \Lambda^2 W=W^\vee$ and the contraction
$V\otimes V^\vee\ra \CZ$ corresponding to the identity of $V$.
Then we have to show that $U$ maps to 0 in $U^\vee\otimes W^\vee$
via the composition
{\small\begin{equation}\label{BtensorB}
U \rxar{g_{U^\vee}}W\otimes V^\vee
\rxar{id_{W}\otimes q}
W\otimes V
\rxar{\otimes B}
(W\otimes W)\otimes (V\otimes V^\vee)\otimes U^\vee
\ra
U^\vee \otimes W^\vee.
\end{equation}}
One sees easily that the above assertion is equivalent to
the property that $B\otimes B$ maps to 0 via the map
$$
(U^\vee\otimes U^\vee)\otimes
(W\otimes W)\otimes
(V^\vee\otimes V^\vee)
\rxar{\ id\otimes p_{W^\vee}\otimes (\neg Q) \ }
(U^\vee\otimes U^\vee)\otimes
(W^\vee).
$$

The above map  factors through
\begin{equation}\label{BtensorB2}
(U^\vee\otimes U^\vee)\otimes
(\Lambda^2 W)\otimes S^2(V^\vee),
\end{equation}
and we have to show that the image of $B\otimes B$
in this space maps to 0 via $ id \otimes id \otimes (\neg Q)$.
\medskip

Write $U^\vee\otimes U^\vee$ as a direct sum
$S^2(U^\vee)\oplus \Lambda^2(U^\vee)$.
By the definition of  $Q$  we get $0$ for
  the contraction $\neg Q$ with $Q$ of the component in
$(\Lambda^2 U^\vee)\otimes (\Lambda^2 W)\otimes S^2(V^\vee)$
of the image of $B\otimes B$.

On the other side, the component of the image of $B\otimes B$ in
$(S^2 U^\vee)\otimes (\Lambda^2 W)\otimes S^2(V^\vee)
=Hom\left(S^2V,(S^2 U^\vee)\otimes (\Lambda^2 W)\right)$
is also zero, because  $S^2(V)$ maps to
$S^2(U^\vee\otimes W)=
\left(\Lambda^2 U^\vee\otimes \Lambda^2 W \right)
\oplus \left( S^2U^\vee \otimes S^2 W \right)$.
\qed

\bigskip
We now want to relate the method to construct such bundles $\E$
as kernels with the direct image method illustrated in section 2.

Consider therefore a tensor
$$\widehat B\in V\otimes {\widehat U}^\vee
\otimes {\widehat W},$$
  and  apply to it the direct image method
of section 2 with twist
$t=2$
(assuming of course that  $\widehat B$ defines
a complete intersection
$\Gamma\subset \PP(V)\times\PP(U')$).
Exact sequence (\ref{Ehat}) gives
\begin{gather}
0 \ra \wedge^2(\widehat W ^\vee) \otimes\cO_{\PP^3}(-2) \ra
\widehat W^\vee\otimes \widehat U^\vee\otimes \cO_{\PP^3}(-1) \ra
\widehat\E_{\widehat B} \ra 0\\
0 \ra \widehat\E_{\widehat B} \ra S^2 (\widehat U^\vee)\otimes
\cO_{\PP^3}\ra \widehat\G_{2\,\widehat B} \ra 0.
\end{gather}
and $\hat\E$ is a vector bundle on $\PP^3$.

Denote
$\widehat \E_{\widehat B}$ simply by $\widehat \E$,
and consider $\hat\E^\vee$: we want to show that there is a tensor $B$
such that $\hat\E^\vee = \E_B.$

Indeed we can dualize the first exact  sequence above,
obtaining
\begin{equation}
0 \ra \hat\E^\vee \ra
{\widehat W}\otimes {\widehat U} \otimes \cO_{\PP'}(1) \ra
{\widehat W}^\vee \otimes\cO_{\PP'}(2) \ra 0.
\end{equation}
Thus $\hat\E^\vee$ is a rank 6 vector bundle and,
by looking at the long   exact cohomology sequences associated
to the twists of the previous exact sequence,
we obtain that the only non-vanishing intermediate cohomology groups of
$\hat\E$ are the two groups
{\small
$$
H^1(\hat\E^\vee(-2))= {\widehat W}^\vee
$$
$$
H^1(\hat\E^\vee(-1))=
\coker\left({\widehat W}\otimes {\widehat U}
\ra {\widehat W}^\vee\otimes V\right)
\cong \left(\ker({\widehat W}\otimes V^\vee \ra
{\widehat W}^\vee\otimes {\widehat U}^\vee)\right)^\vee
$$
}
Thus first of all $\hat\E^\vee$ is again a vector bundle on $\PP^3$ with
Chern polynomial $1+3t+6t^2+4t^3$ and minimal cohomology.

\smallskip

Observe that, in terms of the cross-product involution
applied to the  4-tuple
$$({\widehat W}^\vee, {\widehat U}^\vee,V,\widehat B),$$ we have:
$\left(\ker({\widehat W}\otimes V^\vee \ra
{\widehat W}^\vee\otimes {\widehat U}^\vee)\right)^\vee
= ({\widehat U '})^\vee$.

Hence, if we set
$$U : = {\widehat W}^\vee , W : = ({\widehat U} ' )^{\vee}, B : =
({\widehat B})' , $$
the bundle $\E_B$ associated to $B$ via the kernel construction will be
isomorphic to the bundle $\hat\E^\vee$, as we wanted.

\section{Semistability and moduli space}

In this section we shall show that the explicit geometric
construction we gave before
lends itself to construct a natural moduli space
$\frak A ^0$ for the vector bundles considered in this paper.

Since moduli space for vector
bundles have been constructed in great generality by Maruyama, it
seems natural to
investigate their Gieseker stability (we refer  to \cite{oss} and
especially to \cite{hl}
as  general references).
We conjecture that our bundles are Gieseker stable,
but unfortunately for the time being we only managed to prove
their  slope (Mumford-Takemoto) semistability.

We are however able to prove that our vector bundles are simple,
and we observe then (cf. Theorem 2.1 of \cite{kob})
that moduli spaces of simple vector bundles exist
as (possibly non Hausdorff) complex analytic spaces.

We show indeed that the above moduli space exists as
an algebraic variety.
More precisely, we show that, under a suitable open condition, we can
construct a G.I.T. quotient $\frak A^0$ which is a coarse moduli space.

\bigskip
Recall lemma \ref{E=beta}:
it will lead to a characterization
of the vector bundles obtained from the
kernel construction as an open set in any family of
vector bundles with the above Chern polynomial.

\begin{prop}\label{sstable}
Consider a  rank 6 vector bundle of  $\E$ with total Chern class
$1 + 3 t + 6 t^2 + 4 t^3$,
such that
\begin{enumerate}
\item
$h^0(\E) = 6$
\item
the 6 sections generate a rank 6 trivial subsheaf
with quotient $\tau$
\item
$h^0(\E^\vee) = 0$
\item
$ \E$ is a subbundle of $ 3 \Omega^1 (2)$.
\end{enumerate}
Then $\E$ is  slope-semistable.
         \end{prop}

\proof

Let $\E''$ be a destabilizing subsheaf of rank $r \leq 5$ and
maximal slope $\mu = d/r$ :
without loss of generality we may assume that $\E''$ is  is a
saturated reflexive subsheaf,
and similarly
$\tilde{\E} : = \E'' \cap 6 \hol$ is a saturated reflexive subsheaf
of $6 \hol$.
$$
\xymatrix{
0 \ar[r] &6\hol \ar[r] &\E \ar[r] &\tau \ar[r] &0\\
0 \ar[r] &\tilde\E \ar[r] \ar[u] &\E'' \ar[r] \ar[u] &\tau'' \ar[r] \ar[u] &0\\
&0 \ar[u] &0 \ar[u]\\
}
$$

The slope $ \mu(\E) $ equals $ 1/2$. On the other hand, by hypothesis 4
and since  $\Omega^1 (2)$ is a stable bundle (cf. 1.2.6 b , page 167 of
\cite{oss}), the slope of $\E''$ is at most $2/3$, and
         $ < 2/3$ unless $\E'' \cong \Omega^1 (2)$.

{\bf CLAIM:} $\E$ contains no subsheaf isomorphic to $\Omega^1 (2)$.

{\em Proof of the claim:} $ h^0 (\Omega^1 (2)) = 6 = h^0 (\E)$,
thus this calculation contradicts hypothesis 2.
\qed

\medskip
We have that $ d : = c_1 (\E'') = c_1 (\tilde{\E}) + c_1 (\tau'')$,
and $\tau'' \subset \tau$ is a coherent subsheaf supported on a divisor, thus,
         $ c_1 (\tau'')\leq   c_1 (\tau) = 3$.

On the other hand, $ c_1 (\tilde{\E}) \leq 0$, and if equality holds,
then  $ \tilde{\E} \cong r \hol$.

Hence, $ 1 \leq d  \leq 3$, and we have
$$  2/3 >  \mu = d/r > 1/2  \Leftrightarrow   4 d > 2 r > 3d . $$

These inequalities leave open only the case $ d=3, r=5$.

We show that this case does not exist.

In fact, otherwise we  consider the  quotient by the subbundle
$ \tilde{\E} \cong r \hol$. By hypothesis 3, and the proof of
lemma \ref{unicity} we see that $ \E / \tilde{\E}$ is
an extension corresponding to a tensor of maximal rank, hence
it yields a vector bundle $\cal V$ (cf. corollary \ref{vb}).

Since the torsion sheaf $\tau'' \subset \cal V$, we obtain
$\tau'' = 0$, hence $ d \leq 0$, absurd.

\QED

\begin{rem}
The possible exceptions to slope-stability, in view of the inequalities
$$  2/3 >  \mu = d/r  \geq 1/2  \Leftrightarrow   4 d  \geq 2 r > 3d  $$
are:

1. $  d=1, r=2 $

2. $  d=2, r=4 $.

Matei Toma pointed out how case 2. could be excluded using
Bogomolov's inequality for stable bundles, as done in Lemma 3.1 of his
paper \cite{tom}. The case $r=2$, $c_1 (\tilde{\E}) = -2$ seems
as of now the most difficult case.

Observe that slope-stability of $\E$  implies Gieseker stability of
$\E$, which in turn
implies that there is a point in the moduli space of Gieseker
semistable bundles
corresponding to the
isomorphism class  of $\E$.
\end{rem}

\begin{lem}\label{simple}
Let $\E$ be a vector bundle as in (\ref{E=ker})
with $h^0(\E)=6$ (equivalently, $h^1(\E)=0$)
        and verifying the second assumption.
Then $hom (\E, \E) = 1$, i.e., $\E$ is simple.
\end{lem}

\proof
We consider the exact sequence
$$0 \ra Hom(\E,6\hol) \ra Hom(\E,\E) \ra Hom(\E,\tau) \ra Ext^1(\E,\hol).$$

We have
$Ext^1(\E,\hol)\cong H^1(\E^\vee)\cong H^2(\E(-4))$
and from the exact sequence (\ref{E=ker}) we infer $ H^2(\E(-4)) = 0$.
Since
$Hom(\E,6\hol)=0$ by proposition 3.2, it follows that $Hom(\E,\E) \cong
Hom(\E,\tau)$.

We compute $hom(\E,\tau)$ by considering the exact sequence
$$0 \ra Hom(\tau,\tau) \ra Hom(\E,\tau) \ra Hom(6\hol,\tau).$$
Indeed $hom(\hol,\tau)=h^0(\tau)=0$ (since $h^0(\E)=6$) and, since $\tau$ is
$\hol_G$-invertible, we have $hom(\tau,\tau)=1$.
\qed

\begin{lem}\label{dim=19}
Let $\E$ be a  simple  vector bundle of rank 6, with Chern
classes $ c_1(\E) = 3, c_2(\E) = 6, c_3(\E) = 4$.
Then the local dimension of the  moduli space
$ \frak M^{s} ( 6; 3,6,4)$ of simple vector bundles
at the point corresponding to
$\E$  is at least $19$.
\end{lem}

\proof
The moduli space of simple vector bundles exists (cf. \cite{kob},
Theorem 2.1) and it is well known
that the local dimension is at least equal to the expected dimension $ h^1 (
\E^\vee
\otimes
\E) - h^2 (
\E^\vee
\otimes
\E)$. On the other hand, $\E$    simple means that $ h^0 ( \E^\vee \otimes \E)
= 1$, hence follows also that $ h^3 ( \E^\vee \otimes \E) = 0$,
since by Serre duality $ h^3 ( \E^\vee \otimes \E) =
         h^0 ( \E^\vee \otimes \E (-4)) = 0$.

Thus the
expected dimension  equals $ - \chi ( \E^\vee \otimes \E)  + 1 $ and there
         remains to calculate $ - \chi ( \E^\vee \otimes \E) $.
This can be easily calculated in the case where we have an exact  sequence
$ 0 \ra \E \ra 9 \hol (1) \ra 3 \hol(2)$. We omit the rest of the
easy calculation.
\qed

\medskip
In the following theorem we shall phrase the geometric
meaning of the cross-product involution in terms of a
birational duality of
moduli space of vector bundles, $\frak A ^0$ on $\PP^3$,
${\frak A ^0_*}$ on ${\PP^3}^\vee$.

\medskip
{\bf Main Theorem} {\em
Consider the moduli space $ \frak M^s ( 6; 3,6,4)$ of
rank 6 simple vector
bundles $\cal E$ on $\Pn 3 : = \Proj (V)$
with Chern polynomial $ 1 + 3t + 6 t^2 + 4 t^3$,
and inside it the  open set  $\ \frak A$  corresponding to the simple
bundles with minimal cohomology, i.e., those with
$$
\begin{array}[c]{llll}
(1) & H^i(\E) = 0 \ \forall i \geq 1;
\qquad & (2) & H^i(\E (-1)) = 0 \ \forall i \neq 1;\cr
(3) & H^i(\E (-2)) = 0 \ \forall i \neq 1;
\qquad & (4) &H^i(\E (-3)) = 0 \ \forall i;\cr
(5) & H^i(\E (-4)) = 0 \ \forall i.
\end{array}
$$
%

Then   $\frak A$  is irreducible of dimension 19
and it is bimeromorphic to $ \frak A^0 $,
where $\frak A^0$ is an open set of the G.I.T. quotient space of the
projective space $\frak B$ of
tensors of type $(3,4,3)$,
$\frak B : = \{B\in \PP ({U}^\vee\otimes W \otimes V^\vee)\}$
by the natural action of $SL(W) \times SL(U)$ (recall that $U,W$ are two
fixed vector spaces of dimension 3, while $ V = H^0( \PP^3, \hol (1))$.

Let  moreover $ [B] \in \frak A ^0$ be a general point:
then to $[B]$
corresponds a vector bundle $\E_B$ on $\PP^3$ 
via the kernel construction, 
and also a vector bundle $\E^*_B$ on ${\PP^3}^\vee$, 
obtained from the direct construction
applied to the tensor $B\in U^\vee \otimes W \otimes V^\vee$
(cf. definition \ref{dir_cos} applied to $B$, or equation (\ref{*})).
$\E^*_B$ is the vector bundle $\E_{B'}$, 
where $B'\in W \otimes {U'}^\vee \otimes V$
is obtained from $B$ via the reversing construction
and $[B'] \in {\frak A ^0_*}$.
}\medskip

\proof
To any such tensor $B$ we tautologically associate two linear maps which
we denote by the same symbol,
$$B :  U \otimes V \ra W , \quad B :  U \otimes V \otimes  \hol(1)\ra W \otimes
\hol(1)$$
and using the Euler sequence we define a coherent sheaf $\E$ on $\PP^3$
as a kernel, exactly as in the exact sequence (\ref{E=ker2}) 
(except that surjectivity holds only for $B$ general), 
following what we called the kernel construction.

As we already saw in (\ref{B=beta}), this is equivalent to giving
$\E$ as the kernel of a homomorphism $\beta$ as in (\ref{E=ker}).
Observe that $GL(W) \times GL(U)$ acts on the vector space of such
tensors, preserving the isomorphism class of the sheaf thus obtained.

We define   $\frak B' $ as the open set in  $\frak B $ where
$\beta$ is surjective (thus $\E$ is a rank 6 bundle) and $h^0(\E)=6$. 
Both conditions amount to the surjectivity
of $h^0(\beta)=h^0(B\oplus \epsilon)$, cf. (\ref{E=ker2}),
and imply that $\E$ is a bundle with minimal cohomology,
in the sense of lemma \ref{E=beta}.
We further define
$\frak B ''$ as the smaller open set where the second assumption
is verified,
and we observe then that  lemma \ref{simple}
ensures the existence of a
morphism
$\frak B '' \ra \frak A$ which factors through the action of
$ SL(W) \times SL(U)$.

\smallskip
Since we want to construct a G.I.T. quotient of an open set of $\frak B$,
we let $\frak B^*$ the open set of tensors $B$ whose determinant
defines a cubic surface $G^* \subset {\PP^3 }^\vee$, i.e., we have an
exact sequence on ${\PP^3 }^\vee$ of the form (set $\hol_* : =
\hol_{{\PP^3 }^\vee}$)
\begin{equation}\label{*}
0 \rar U\otimes \cal O_*(-1)  \rxar{B} W\otimes \cal O_*  \rar
\cal G^* \rar 0.
\end{equation}
Since the determinant map is obviously $SL(W) \times SL(U)$-invariant,
the tensors in  $\frak B^*$ are automatically semistable points
for the $SL(W) \times SL(U)$-action, by virtue of the criterion of
Hilbert-Mumford.

Observe now that the maximal torus $\C^* \times \C^*$ of
$ GL(W) \times GL(U)$ acts trivially on $\frak B$, thus
we get an effective action of $SL(W) \times SL(U)$ only upon
dividing by a finite group $  K' \cong (\Z/3)^2$.

We claim that  $(SL(W) \times SL(U)) / K'$ acts freely on
the open subset $\frak B^{**}   \subset \frak B^* $,
$\frak B^{**} = \{ B \in\frak B^*  | End(\cal G^* ) = \C \}$.

This is clear since the stabilizer of $B$ corresponds uniquely
to the group of automorphisms of $ \cal G^* $,
and any such automorphism acts on $W \cong H^0(\cal G^* )$,
and induces a unique automorphism of $U$ in view of the exact
sequence (\ref{*}).  But every automorphism is multiplication
by a constant, thus it yields an element in $K'$.

We want to show that the orbits are closed.
But the orbits  are contained in the
fibres of the determinant map: thus,
it suffices to show that, fixed the cubic surface $G^*$, if
we have a 1-parameter family where $ \cal G_t \cong \  \cal G_1$
for $ t \neq 0$, then also $ \cal G_0 \cong \  \cal G_1$.

This holds on the smaller open set $\frak B^{***} \subset \frak B^{**}$
consisting of the tensors such that the cubic surface $G^*$
is smooth: since then $ \cal G_0$ is invertible, and the Picard group
of  $G^*$ is discrete.

We have proven that $\frak B^{***}$ consists of stable points,
and observe that  the condition $End(\cal G^* ) = \C $
holds if $\cal G^* $ is $\hol_{G^*}$-invertible, or
it is torsion free and $ G^*$ is normal.
Therefore the open set $\frak B^{st}$ of stable points is nonempty.

\smallskip
We define $\frak A ^0$  as the open set of the G.I.T. quotient
corresponding to
$\frak B^{st} \cap \frak B ''$.

\medskip
The fact that $\frak A$ is irreducible follows since
every bundle $\E$ in $\frak A$ has a cohomology table which
(by Beilinson's theorem, as explained in lemma \ref{E=beta})
implies that $\E$ is obtained from a tensor
$B$ in the open subset ${\frak B '}^0 \subset {\frak B '}$
consisting of those $B$ for which the corresponding bundle 
$\E$ is simple
(note that ${\frak B '}^0  \supset {\frak B '' }$).

Now,  $\dim \frak A^0 = 19 $, while $\dim \frak A \geq 19 $
by \ref{dim=19}; we only need to observe that if
$[B] , [B'] \in \frak A^0$ and two bundles
$\E_B$ and $\E_{B'}$ are isomorphic, then the corresponding
tensors  $B, B'$ are $ GL(U) \times GL(W)$ equivalent, since
they express the multiplication matrix for the intermediate cohomology
module $ H^1_*( \E)$. Thus $[B] = [B'] \in \frak A^0$.

It follows on the one side that
$\frak A^0 $  parametrizes isomorphism classes of bundles,
and on the other side  that $\frak A^0 $
maps bijectively to an open set in $\frak A$, in particular
$\dim \frak A = 19 $, since $\frak A $ is irreducible.

\QED
\bigskip

{\bf Acknowledgement}
The  authors would like to thank I. Dolgachev for  pointing out
how the Schur quadric could be related to the cross-product involution, 
for suggesting a third construction of the bundles $\E$, 
and for  hospitality and the nice atmosphere during  the visit 
of the first author in Ann Arbor in january 2006. 
The present research was carried out
in the realm  of the DFG Schwerpunkt ``Globale Methoden in der
komplexen Geometrie''. The second author also profited from a travel
grant  from a DAAD-VIGONI program.



\end{document}